\documentclass[11pt]{article}          
\usepackage{geometry}
\usepackage{hyperref}
\usepackage{algorithm,algorithmic}
\usepackage{graphicx}
\usepackage{amsmath}
\usepackage{amsfonts}
\usepackage{latexsym}
\usepackage{comment}
\usepackage{color}
\usepackage{listings}
\usepackage{placeins}
\usepackage{subcaption}
\usepackage[labelfont=bf]{caption}  
\usepackage{multirow}
\ifpdf
  \DeclareGraphicsExtensions{.eps,.pdf,.png,.jpg}
\else
  \DeclareGraphicsExtensions{.eps}
\fi

\usepackage{epstopdf}
\usepackage{enumitem}

\usepackage{amsopn}

\newcommand {\vv}  { {\bf v} }

\newcommand{\hf}{\frac12}

\newcommand{\bfv}{ {\bf{v}}}

\newcommand{\bbR}{\mathbb{R}}
\newcommand{\bbC}{\mathbb{C}}

\newcommand{\bfI}{{\bf I}}

\newcommand{\bfQ}{{\bf Q}}

\newcommand{\bfd}{{\bf d}}
\newcommand{\bfe}{{\bf e}}

\newcommand{\bfr}{{\bf r}}

\newcommand{\bfy}{{\bf  y}}
\newcommand{\bfx}{{\bf  x}}

\newcommand{\bfz}{{\bf z}}

\newcommand{\Fcurly}{\mathcal{F}}

\newcommand{\bfepsilon}{\boldsymbol \epsilon}

\usepackage{mathtools}
\usepackage{url}

\usepackage{array}

\usepackage{color}

  \newcommand{\doneone}[2][]{}

\newcommand{\removed}[1]{}
\newcommand{\done}[2][]{}
\newcommand{\mdone}[1]{}

\newcommand{\ds}{{\displaystyle}}

\def \IhH         {\bfI_h^H}

\def \IHh         {\bfI_H^h}

\newlength{\myrowheight}
\date{\today}

\begin{document}

\title{Multigrid Optimization for Large-Scale Ptychographic Phase Retrieval} 



\author{%
Samy Wu Fung\thanks{Department of Mathematics, Emory University,  Atlanta, GA. \texttt{samy.wu@emory.edu}, \newline \url{https://sites.google.com/site/samywufung/}}
\and
Zichao (Wendy) Di \thanks{Mathematics and Computer Science Division, Argonne National Laboratory, Lemont, IL. \texttt{wendydi@anl.gov}, \newline \url{http://www.mcs.anl.gov/person/zichao-wendy-di}}
}







\maketitle

\begin{abstract}
Ptychography is a popular imaging technique that combines diffractive imaging with scanning microscopy. 
The technique consists of a coherent beam that is scanned across an object in a series of overlapping positions, leading to reliable and improved reconstructions.
Ptychographic microscopes allow for large fields to be imaged at high resolution at the cost of additional computational expense.
In this work, we propose a multigrid-based optimization framework to reduce the computational burdens of large-scale ptychographic phase retrieval. 
Our proposed method exploits the inherent hierarchical structures in ptychography through tailored restriction and prolongation operators for the object and data domains.
Our numerical results show that our proposed scheme accelerates the convergence of its underlying solver and outperforms the Ptychographic Iterative Engine (PIE), a workhorse in the optics community.
\vspace{2mm}

\end{abstract}

\section{Introduction}
\begin{figure}[t]
    \centering
    \includegraphics[width = 0.7\textwidth]{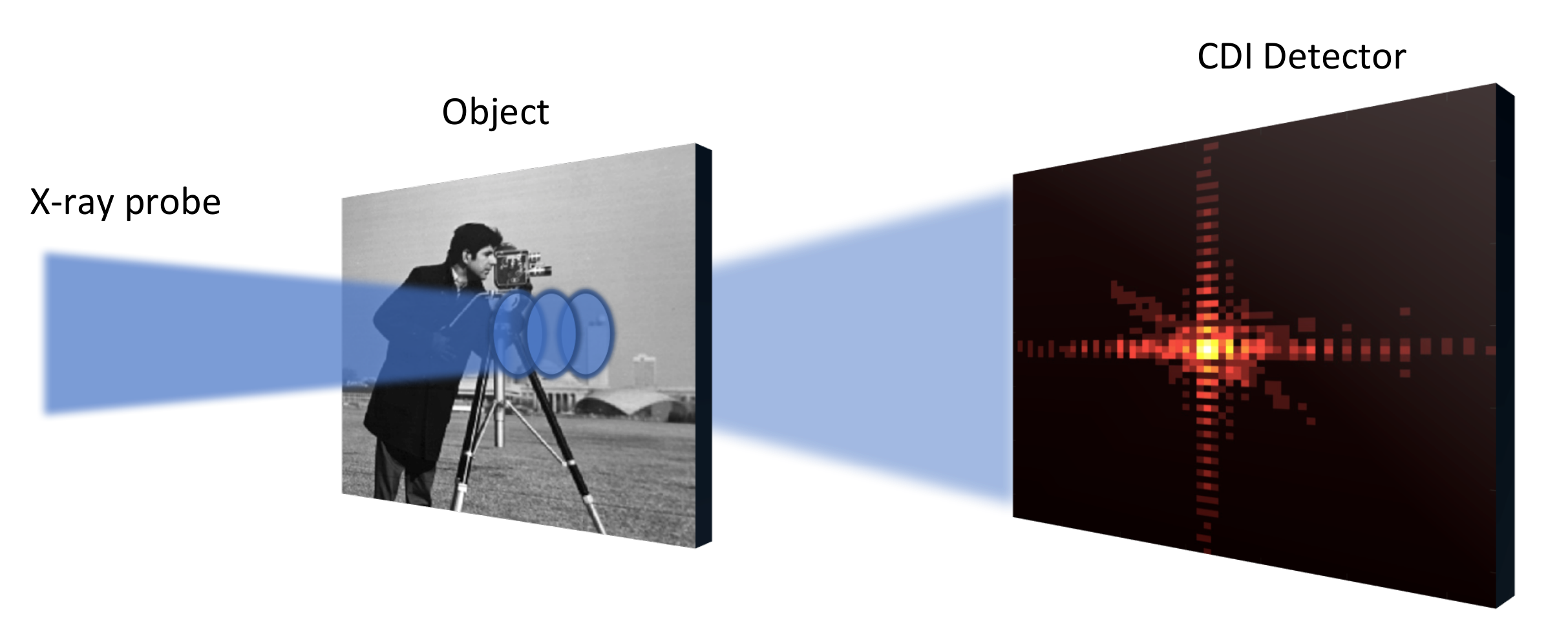}
    \caption{\textit{Schematic illustration of the ptychography experiment with three successive overlapping scans.}}
    \label{fig:ptychoExperiment}
\end{figure}

Ptychography is a coherent diffractive imaging (CDI) technique that arises in applications such as materials science \cite{holler2017high,hoppe2013high,pelz2014fly}, biology \cite{marrison2013ptychography,suzuki2016dark}, and x-ray crystallography \cite{de2016ptychographic}. The technique was originally proposed to improve the resolution in electron or x-ray microscopy by replacing single-element detectors with two-dimensional array detectors, combining diffractive imaging with scanning microscopy \cite{gursoy2017direct}. More precisely, a coherent beam is scanned across an object in a series of overlapping positions that couples information between successively collected diffraction patterns (see Fig.~\ref{fig:ptychoExperiment}). As is well known in single-pattern CDI, finite support constraints are crucial for the convergence of classical algorithms such as error reduction \cite{gerchberg1972practical}, hybrid input output (HIO) \cite{fienup1982phase}, gradient-based algorithms \cite{candes2015phase}, relaxed averaged alternating reflections \cite{luke2004relaxed}, and saddle point optimization \cite{marchesini2007phase,tripathi2015visualizing}. In ptychography, \textit{a priori} knowledge of the scanning positions automatically delivers these constraints, leading to faster and more robust reconstructions than single-pattern CDI \cite{bunk2008influence,deng2015continuous, huang2014optimization, maiden2009improved}.


%

The widespread application of ptychography has led to considerable research on methods for its numerical reconstruction. Among the most popular methods is the ptychographic iterative engine (PIE) \cite{maiden2009improved,rodenburg2008ptychography,rodenburg2004phase}. The algorithm consists of alternating projections onto non-convex modulus constraint sets and is popular in the optics community \cite{de2016ptychographic, giewekemeyer2011ptychographic, marchesini2007invited,marrison2013ptychography,suzuki2016dark}
Mathematically, PIE is equivalent to a projected steepest descent algorithm applied to a particular error metric (see Sec.~\ref{sec:mathFormulation}).
Other well-known approaches include standard gradient-based techniques \cite{qian2014efficient, yang2011iterative}, and Wirtinger flow \cite{candes2015phase,xu2018accelerated}, which uses careful initialization via a spectral method. Recently, a lifting approach (PhaseLift) was introduced, where the phase retrieval problem is reformulated as a convex optimization problem at the expense of solving for a quadratically increased number of unknown variables \cite{candes2015phaseReview}; unfortunately, this quickly becomes intractable for large-scale problems. Despite serious efforts, large-scale ptychographic phase retrieval continues to be challenging as high resolution demands and small scanning beams generate large volumes of data. 

In this work, we consider the multigrid-based optimization framework (MG/OPT) presented in \cite{nash2000multigrid} for solving the ptychographic phase retrieval. The MG/OPT scheme is a general-purpose framework designed to accelerate the large-scale nonlinear optimization problems. MG/OPT is inspired by the full approximation scheme (FAS) \cite{brandt2011multigrid,trottenberg2000multigrid}, and has stronger convergence guarantees than traditional nonlinear multigrid methods \cite{nash2010convergence}. Moreover, MG/OPT can be extended to optimization problems with equality and inequality constraints \cite{lewis2005model}. When applied to unconstrained optimization problems, MG/OPT is equivalent to applying FAS to the first-order optimality conditions. However, MG/OPT ensures a descent direction in the coarse-grid correction under mild assumptions; this leads to a globally convergent algorithm when a linesearch step is performed \cite{nash2010convergence}.

Our work is motivated by the successes of MG/OPT on different nonlinear problems \cite{Di2012truncated,lewis2005model,nash2000multigrid}, as well as a previously applied multilevel scheme
\cite{seifert2006multilevel}. The multilevel scheme proposed in \cite{seifert2006multilevel} uses a \textit{single slash} cycle on 1D single-pattern CDI phase retrieval problems; specifically, the coarse grid calculation only contributes to a better initial guess for the fine-grid problem.
Our approach differs from \cite{seifert2006multilevel} in that we use multiple V-cycles as well as a linesearch scheme that guarantees a descent direction during the coarse-grid correction. Our work is also motivated by the different hierarchical structures that ptychography exhibits, and which can be exploited by a multilevel scheme as we will show in Sec.~\ref{subsec:MGOPTforPtycho}.

The remainder of this paper is organized as follows. In Sec.~\ref{sec:mathFormulation}, we review the mathematical formulation and some common algorithms for the ptychographic phase retrieval. In Sec.~\ref{sec:MGOpt}, we give an overview of the MG/OPT framework. In Sec.~\ref{sec:numResults}, we illustrate the potential of MG/OPT for solving ptychographic phase retrieval on a few synthetic experiments, and finally, we conclude with a discussion in Sec.~\ref{sec:conclusion}.

\section{Mathematical Background}
\label{sec:mathFormulation}
In this section, we give an overview of the 2D ptychographic phase retrieval as well as some of the popular methods used for its reconstruction.
Since we follow the discretize-then-optimize approach, we limit the discussion to the discrete setting - the continuous formulation of the general phase retrieval problem can be found in \cite{luke2002optical}. 

Let $\bfz=\bfx+\bfy i \in \bbC^{n^2}$ be the object of interest, and $\bfd_k \in \bbR^{m^2}$ be the observed data (or intensities) measured from the $k^{th}$ probe, where $n^2$ and $m^2$ are the dimensions of the vectorized object and data resolution images, respectively. A ptychography experiment is modeled by 
\begin{align} \label{eq:forwardProblem}
    \bfd_k = |\mathcal{F}(\bfQ_k\bfz)|^2 + \bfepsilon_k, \quad k=1,\ldots,N,
\end{align}
where $N$ is the total number of probes (or scanning positions), $\mathcal{F}\colon \mathbb{C}^{n^2} \mapsto \mathbb{C}^{m^2}$ is the two-dimensional discrete Fourier operator, $\bfQ_k \in \bbR^{n^2 \times n^2}$ is the $k^{th}$ probe (a diagonal \textit{illumination} matrix), and $\bfepsilon_k \in \bbR^{n^2}$ is the noise corresponding to the $k^{th}$ measurement error. The diagonal elements of $\bfQ_k$ are nonzero in the columns corresponding to the pixels being illuminated in the object at scanning step $k$. 

There are different ways for formulating the reconstruction problem. One such formulation is the intensity Gaussian error metric \cite{qian2014efficient}
\begin{align}\label{eq:misfit1}
  \min_\bfz \Phi_\mathcal{IG}(\bfz) = \hf \sum_{k=1}^N \| |\mathcal{F}(\bfQ_k \bfz)|^2 - \bfd_k\|_2^2,
\end{align}
which involves minimizing the data misfit. Here, $\Phi_\mathcal{IG}\colon \mathbb{C}^{n^2} \mapsto \bbR$ is a real-valued cost function defined on the complex domain, and is therefore not complex-differentiable \cite{remmert2012theory}. To overcome the lack of complex-differentiability, it is common to employ the notion of $\bbC\bbR$ (Wirtinger) Calculus, where the derivatives of the real and imaginary parts of $\bfz$ are computed independently \cite{remmert2012theory,sorber2012unconstrained}. For these real-valued functions, the mere existence of these Wirtinger derivatives is necessary and sufficient for the existence of a stationary point \cite{brandwood1983complex,remmert2012theory,sorber2012unconstrained}. Using Wirtinger calculus, the partial gradients for~\eqref{eq:misfit1} can be computed as
\begin{equation}
    \begin{split}
        \nabla_\bfx \Phi_\mathcal{IG} = \sum_{k=1}^N \bfQ_k^\top \mathcal{F}\left(\overline{\mathcal{F} \bfQ_k \bfz \odot \bfr_k}\right),
        \quad 
        \nabla_\bfy \Phi_\mathcal{IG} = \sum_{k=1}^N i\bfQ_k^\top \mathcal{F} \left(\overline{\mathcal{F} \bfQ_k \bfz \odot \bfr_k}\right),
    \end{split}
\end{equation}
where $\bfr_k = |\mathcal{F} \bfQ_k \bfz|^2 - \bfd_k$ is the residual of the $k^{th}$ probe , and the resulting Wirtinger derivative is 
\begin{align}
\nabla_\bfz \Phi_\mathcal{IG} = [\nabla_\bfx \Phi_\mathcal{IG}^\top \quad \nabla_\bfy \Phi_\mathcal{IG}^\top]^\top. 
\end{align} 
Popular methods that solve this problem include classical gradient-based algorithms \cite{qian2014efficient,yang2011iterative} and Wirtinger flow \cite{candes2015phase,xu2018accelerated}. Variants of the intensity Gaussian error metric include the amplitude Gaussian metric \cite{chang2018variational,qian2014efficient,wen2012alternating}, intensity Poisson metric \cite{chang2018total,chang2018variational,thibault2012maximum}, and the weighted intensity Gaussian metric \cite{qian2014efficient}, all of which measure some variants of the misfit between the forward model and the observed data.

Another formulation of the inverse problem involves solving a feasibility problem. In particular, let the $k^{th}$ measurement constraint set and its corresponding projection operator be denoted by
\begin{align}
  \mathcal{M}_k = \{\bfz \in \bbC^{n^2} : |\mathcal{F} \left(\bfQ_k \bfz \right)| = \bfd_k \}, \;\; \text{ and } \;\; \mathcal{P}_{\mathcal{M}_k} (\bfz) = \mathcal{F}^{-1} \left[ \sqrt{\bfd_k} \odot \exp \big(i \: \theta(\mathcal{F} \left(\bfQ_k \bfz \right)\big) \right],
\end{align}
respectively, where $\odot$ is the Hadamard (or element-wise) product, the square root in $\mathcal{P}_{\mathcal{M}_k}$ is applied element-wise, and $\theta\colon\bbC \to [-\pi, \pi)$ is the element-wise principal argument function \cite{remmert2012theory}. The feasibility problem is then formulated as
\begin{align}\label{eq:feasibility}
  \text{ find } \bfz \in \mathbb{C}^{n^2} \;\; \text{ such that } \;\;  \bfz \in \bigcap_{k=1}^N \mathcal{M}_k.
\end{align}
Perhaps the most well-known numerical technique for solving~\eqref{eq:feasibility} is the alternating projection algorithm, PIE, and its variants \cite{konijnenberg2016combining,maiden2009improved}. PIE consists of performing projections onto the sets $\mathcal{M}_k$ that are calculated one at a time in a sequential manner, where $\gamma \in \bbR$ is a relaxation scalar that helps suppress the noise effect in the data (see Alg~\ref{alg:PIE}). These methods have enjoyed great success in the optics community \cite{de2016ptychographic,giewekemeyer2011ptychographic,konijnenberg2016combining,marrison2013ptychography,rodenburg2008ptychography}; however, convergence of these methods is not always guaranteed because the sets $\mathcal{M}_k$ are non-convex \cite{luke2002optical}.
\begin{algorithm}[t]
  \begin{itemize}
    \item initialize $N$ = number of probes
    \item for $j=1,2,\ldots$
    \begin{itemize}
      \item for $k=1,2,\ldots,N$

           \hspace{3mm} $\bfz^{(j)} = \bfz^{(j)} + \left(\dfrac{|\bfQ_k|}{\max{(\bfQ_k})} \dfrac{\bfQ_k^T}{|\bfQ_k|^2 + \gamma}\right) \big( \mathcal{P}_{\mathcal{M}_k} (\bfz^{(j)}) - \bfQ_k \bfz^{(j)}\big)$
      \item end
    \end{itemize}
    \item check convergence criteria
  \end{itemize}
 \caption{Ptychographic Iterative Engine (PIE)}
 \label{alg:PIE}
\end{algorithm}

Indeed, problem~\eqref{eq:feasibility} can only be solved in the ideal, noise-free case. As soon as there is noise, we are left to minimize some distance metric based on~\eqref{eq:feasibility}. In fact, when $\gamma=0$ and the probes are binary in Alg.~\ref{alg:PIE}, each update in PIE corresponds to a projected steepest descent iteration that solves 
\begin{align}\label{eq:distanceMisfit}
  \min_\bfz \Phi_\mathcal{M}(\bfz) = \hf \sum_{k=1}^N \| \mathcal{P}_{\mathcal{M}_k} (\bfz) - \bfQ_k \bfz \|_2^2,
\end{align}
where $\Phi_{\mathcal{M}}\colon \bbC^{n^2} \mapsto \bbR$ describes the distance from the current point to the set $\bigcap_{k=1}^N \mathcal{M}_k$. Its corresponding Wirtinger derivative can be shown to be
\begin{align}\label{eq:distanceGradient}
  \nabla_\bfz \Phi_\mathcal{M}(\bfz) = \sum_{k=1}^N \nabla_\bfz \Phi_{\mathcal{M}_k}(\bfz) = \sum_{k=1}^N \bfQ_k  \big( \mathcal{P}_{\mathcal{M}_k} (\bfz) - \bfQ_k \bfz \big)
\end{align}
\cite{barakat1985algorithms,luke2002optical,maiden2017further}. In particular, each inner iteration of PIE in Alg.~\ref{alg:PIE} is equivalent to 
\begin{align}
  \bfz^{(j)} = \bfz^{(j)} - \nabla_\bfz \Phi_{\mathcal{M}_k}(\bfz^{(j)}).
\end{align}
To the best of our knowledge, gradient-based algorithms for ptychography have been mostly used to solve the intensity Gaussian metric~\eqref{eq:misfit1} \cite{candes2015phase,qian2014efficient,xu2018accelerated,yang2011iterative}.
However, it has been observed that optimizing over $\Phi_\mathcal{IG}$ (as shown in Sec.~\ref{subsec:compareMisfits}) is prone to getting stuck in poor-quality local minima, i.e., artifacts in the reconstructions \cite{qian2014efficient,yeh2015experimental}. This has resulted in different variants of $\Phi_{\mathcal{IG}}$ \cite{qian2014efficient,yang2011iterative}, as well as careful spectral initialization using Wirtinger Flow \cite{candes2015phase} to be considered in the general phase retrieval community. In this work, we instead focus on the distance metric $\Phi_\mathcal{M}$ as the objective function for its stability. 

Despite the aforementioned efforts, the ptychographic phase retrieval continues to be computationally demanding, especially when extended to 3D applications \cite{gilles20183d,gursoy2017direct}, where the image sizes grow dramatically.
The need to solve these large-scale problems quickly and accurately as a result of the increasing capacities of ptychographic microscopes thus motivates our pursuit of multigrid techniques for the efficient computations of their solutions. 

\section{Multigrid Optimization}
\label{sec:MGOpt}

MG/OPT is a multigrid-based optimization framework designed to solve large-scale nonlinear optimization problems \cite{nash2000multigrid}. Its goal is to accelerate the convergence of traditional iterative algorithms by exploiting the hierarchy of the original optimization problems. The scheme has shown success in a broad class of problems such as PDE-constrained optimal control problems \cite{lewis2005model} as well as in the generation of centroidal Voronoi tessellations \cite{Di2012truncated}. In this section, we review the MG/OPT framework by introducing some notation, followed by a description of an MG/OPT cycle, and finally with a discussion on its extension to ptychography.

\subsection{Notation} We employ the standard notation used in the multigrid community, where the subscript $h$ denotes the fine grid, and $H$ denotes the coarse grid.
In order to transfer information between grid levels, we denote the restriction operator by $\IhH \in \bbR^{n_H^2 \times n_h^2}$, which projects the variable and the gradient from the fine grid to the coarse grid, and the prolongation operator by $\bfI_H^h \in \bbR^{n_h^2 \times n_H^2}$, which interpolates the search direction from the coarse grid to the fine grid. 
We also denote the underlying optimization algorithm by "OPT". It is assumed to be convergent in the sense that, under appropriate assumptions on the objective function $\Phi$,
\begin{align}
\ds \lim_{ j \rightarrow \infty} \left| {\nabla \Phi(\bfz^{(j)})} \right| = 0,
\end{align}
where $\{\bfz^{(j)}\}$ are the iterates computed by OPT. 

Moreover, we define OPT as a function of the form 
$
\bfz^+ \leftarrow \mbox{OPT}(\Phi,\vv,\bar {\bfz},K)
$
which applies $K$ iterations of OPT to the problem 
\begin{align}
  \ds \min_{\bfz_H} \Phi(\bfz_H) - \vv^{T} \bfz_H 
\end{align}
with initial guess $\bar {\bfz}$ to obtain $\bfz^{+}$. If OPT is required to proceed until convergence, we rewrite this as $ \bfz^+ \leftarrow \mbox{OPT}(\Phi ,v,\bar {\bfz})$.

\subsection{MG/OPT Cycle} \label{subsec:MG/OPTcycle}
Given an initial estimate of the solution ${\bfz}_h^{(0)}$ on the fine grid, set $\bfv_h = 0$.  Select non-negative integers $k_1$ and $k_2$ satisfying $k_1+k_2>0$.  Then for $j = 0, 1, \ldots$, set
\begin{align}
\bfz_h^{(j+1)} \leftarrow
\mbox{MG/OPT}(\Phi_h,\bfv_h,\bfz_h^{(j)}),
\end{align}
where the function MG/OPT is defined as follows.
\begin{itemize}             
  \item {\em Coarse-grid solve:}  If on the coarsest grid, then solve the optimization problem:
   \[ 
    \bfz_h^{(j+1)} \leftarrow
    \mbox{OPT}(\Phi_h,\bfv_h,\bfz_h^{(j)}).
    \]
    Otherwise,
     \setlength\itemsep{0.4em}
  \item {\em Pre-smoothing:}
    \[
    \bar {\bfz}_h \leftarrow \mbox{OPT}(\Phi_h,\bfv_h,\bfz_h^{(j)},k_1)
    \]
     \setlength\itemsep{0.4em}
  \item {\em Coarse-grid correction:}
     \setlength\itemsep{0.4em}
    \begin{itemize}
      \item Compute
        \begin{eqnarray*}
          \bar {\bfz}_H &=&\IhH \bar{\bfz}_h  \\
          \bar \vv & = & \IhH \bfv_h +\nabla \Phi_H (\bar {\bfz}_H)-\IhH \nabla \Phi_h(\bar {\bfz}_h)
        \end{eqnarray*}
         \setlength\itemsep{0.4em}
      \item Apply MG/OPT recursively to the surrogate model:
        \[
        \bfz_H^{+} \leftarrow
        \mbox{MG/OPT}(\Phi_H,\bar  \vv,\bar \bfz_H)
        \]
         \setlength\itemsep{0.4em}
      \item Compute the search directions \[\bfe_H = \bfz_H^{+}  - \bar \bfz_H\] and \[\bfe_h = \IHh \bfe_H.\]
         \setlength\itemsep{0.4em}
      \item Use a linesearch to determine \[\bfz_h^{+} ={\bar \bfz}_h+ \alpha \bfe_h\] satisfying $\Phi_h(\bfz_h^{+}) \le \Phi_h({\bar \bfz}_h)$.
    \end{itemize}
     \setlength\itemsep{0.4em}
  \item {\em Post-smoothing:}
    \[
    \bfz_h^{(j+1)} \leftarrow \mbox{OPT}(\Phi_h,\bfv_h,\bfz_h^+,k_2)
    \]
\end{itemize}

\subsection{Extension to Ptychography}
\label{subsec:MGOPTforPtycho}

We now describe MG/OPT in the context of ptychography for the distance error metric~\eqref{eq:distanceMisfit}. For brevity, we consider a two-grid hierarchy, however, the discussion presented below can naturally be extended to multiple levels of grids as well as the 3D case \cite{gilles20183d}. 

In ptychography, the fine-grid problem in MG/OPT is given by 
\begin{align}
  \min_{\bfz_h} \Phi_{\mathcal{M},h} (\bfz_h) = \hf \sum_{k=1}^N \left\| \mathcal{P}_{\mathcal{M}_{k,h}} (\bfz_h) - \bfQ_{k,h} \bfz_h \right\| _2^2
  ,
\end{align}
where 
\begin{align}
  \mathcal{P}_{\mathcal{M}_{k,h}} (\bfz_h) = \Fcurly^{-1} \left( \sqrt{\bfd_{k,h}} \odot \exp \big(i \: \theta\left(\Fcurly\left(\bfQ_{k,h} \bfz_h\right)\right)\big) \right),
\end{align} 
and the additional subscript $h$ represents a certain fine grid level. The coarse-grid surrogate problem can then be written as 
\begin{equation}
\begin{split}
  \min\limits_{\bfz_H} \;\; & \Phi_{\mathcal{M},H}(\bfz_H) - \vv^{T} \bfz_H  \\
  \vspace{0.5mm}
  = & \hf \sum\limits_{k=1}^N \left\| \mathcal{P}_{\mathcal{M}_{k,H}} (\bfz_H) - \bfQ_{k,H} \bfz_H \right\| _2^2- \left( \nabla_{\bfz_H} \Phi_{\mathcal{M},H}(\IhH \bar{\bfz}_h)-\IhH \nabla_{\bfz_h} \Phi_{\mathcal{M},h}(\bar{\bfz}_h)\right)^T \bfz_H, \\
\end{split}
\end{equation} 
where 
\begin{align}
\mathcal{P}_{\mathcal{M}_{k,H}} (\bfz_H) = \Fcurly^{-1} \left( \sqrt{\tilde{\bfI}_h^H \bfd_{k,h}} \exp \Big(i \: \theta\left(\Fcurly\left(\IhH\bfQ_{k,h} \bfz_H\right)\right)\Big) \right).
\end{align}
Here, we also add a subscript $h$ or $H$ to the grid-independent objective function $\Phi_\mathcal{M}$ to distinguish the objectives accross grid levels.

\begin{figure}[t]
  \centering
  \begin{tabular}{ccc}
  \textbf{a)} low-pass filter scheme
  &
  \textbf{b)} original reconstruction
  &
  \textbf{c)} cropped reconstruction
  \\
  \includegraphics[width = 0.355\textwidth, height = 1.5in]{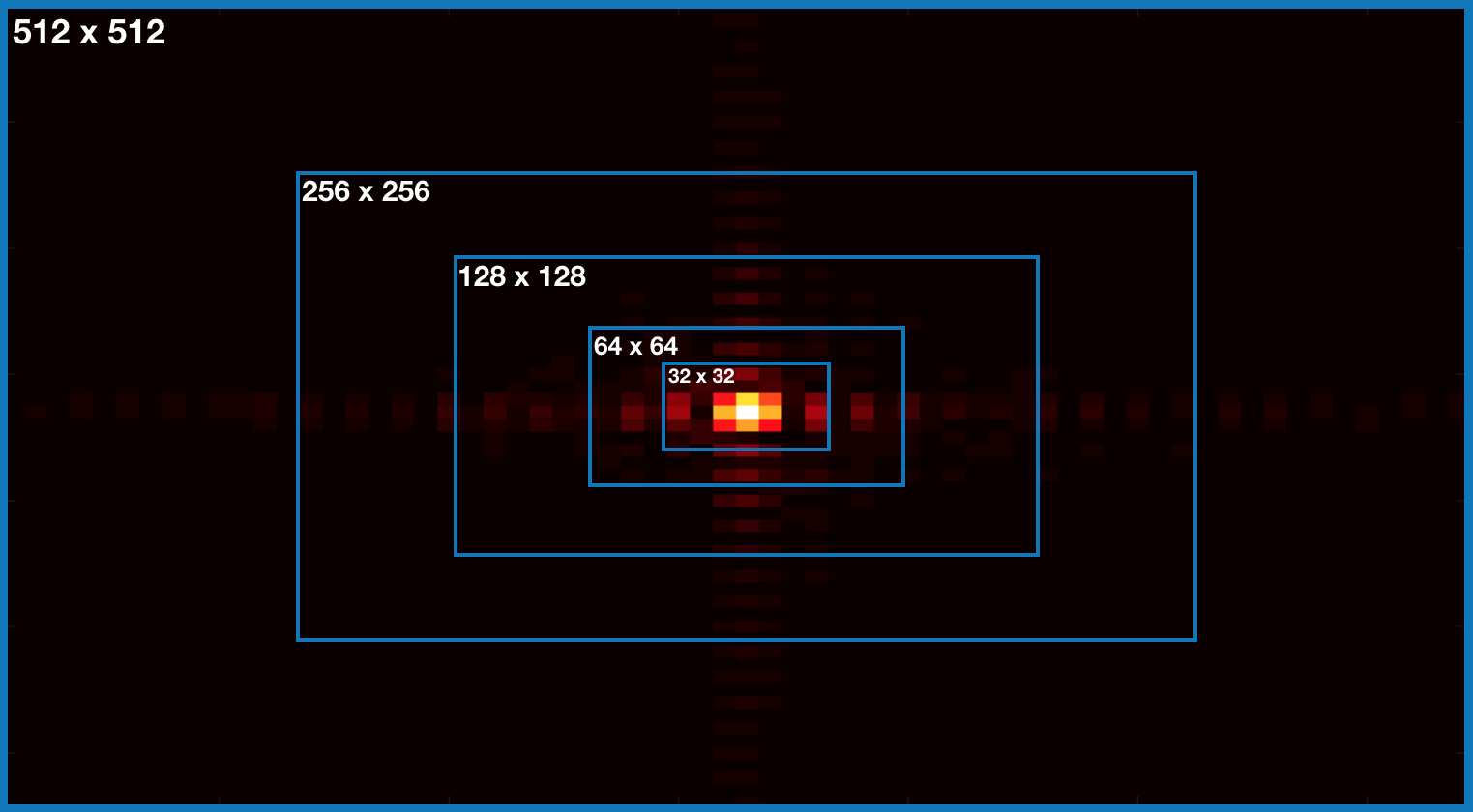}
  &
  \includegraphics[width=0.25\textwidth, height=1.5in]{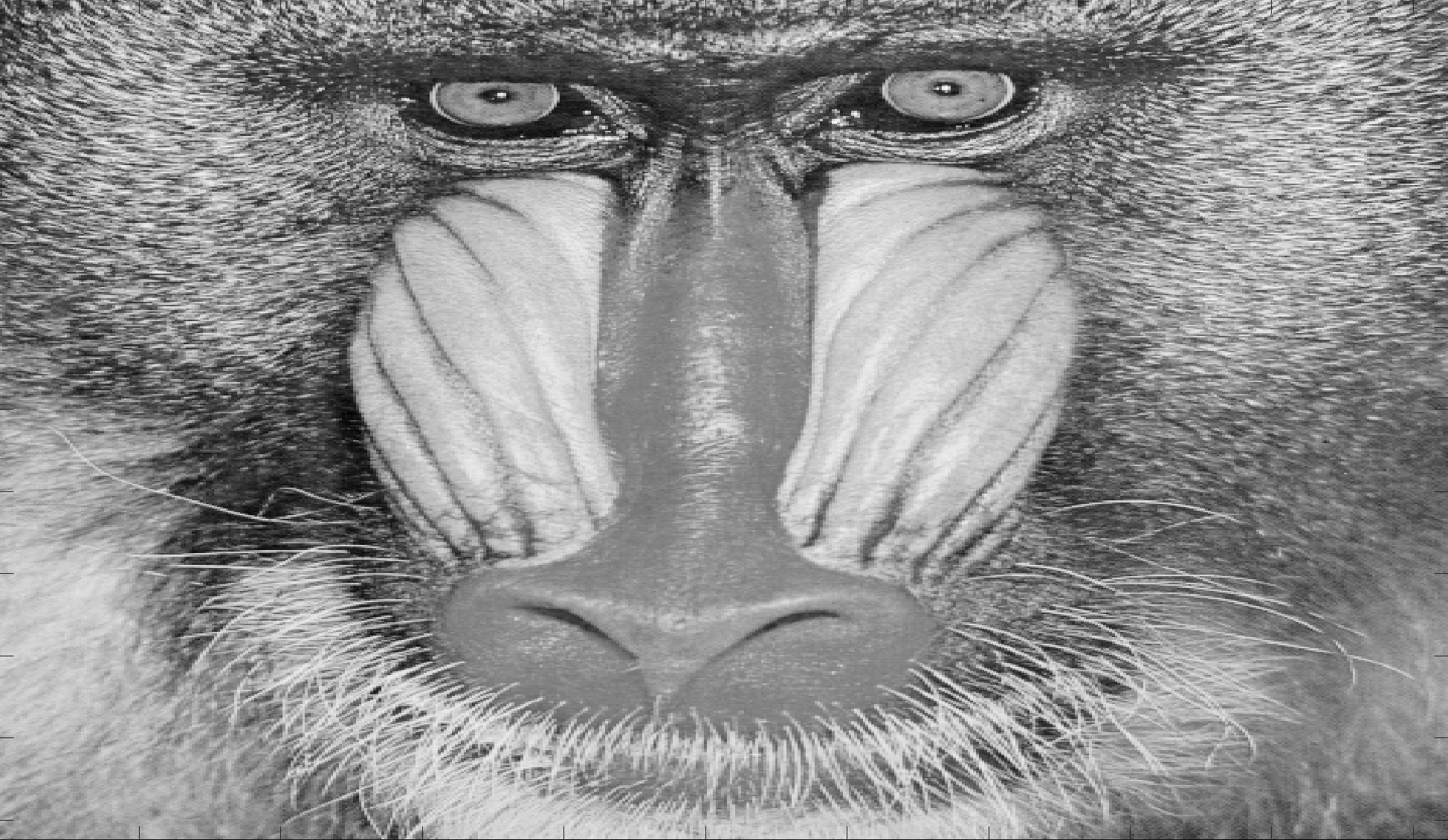}
  &
  \includegraphics[width=0.25\textwidth, height=1.5in]{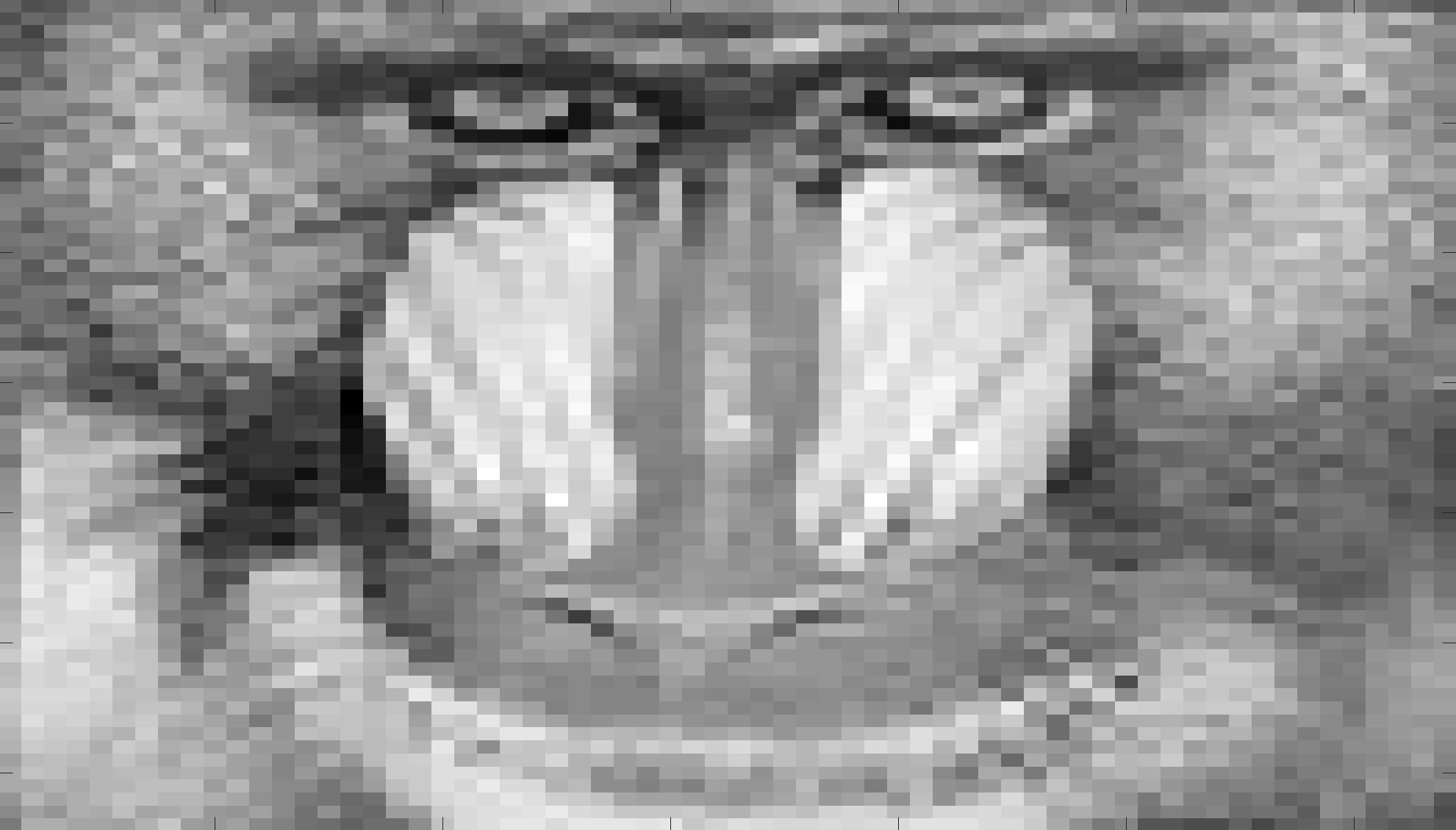}
  \end{tabular}
  \caption{\textit{\textbf{a)} Data coarsening at different grids. 
  \textbf{b)} Reconstruction of linear problem from original ($512 \times 512$) diffraction pattern. \textbf{c)} Reconstruction of linear problem from cropped high-frequencies ($64 \times 64$) diffraction pattern.}}
  \label{fig:subsampledData}
  \vspace{-5mm}
\end{figure}

Unlike the standard MG/OPT where only the variable and gradients are transversed across levels, we also consider transversing the data via the data-restriction operator $\tilde{\bfI}_h^H\in \bbR^{m_H^2 \times m_h^2}$. 
To maintain the data fidelity across different grids, we employ the \textit{low-pass filtering} concept \cite{antoniou2016digital,gonzalez2002digital}, where we only sample the low-frequency components of the data for the coarse-grid problem. More precisely, the data-restriction operator $\tilde{\bfI}_h^H\in \bbR^{m_H^2 \times n_h^2}$ is a binary diagonal matrix 
\begin{align}
  [\tilde{\bfI}_h^H]_{j,j}= 
  \begin{cases}
    1 & \text{if $j$ corresponds to the chosen low-frequency pixel} \\
    0 & \text{otherwise},\\
  \end{cases}
\end{align}
which crops out the high-frequency components and keeps the low-frequency components; see Fig.~\ref{fig:subsampledData}a. This choice of restriction is motivated by the fact that in the frequency domain, the high frequency modes correspond to the detailed features of the object (e.g., the sharp edges) and can only be captured by fine resolutions on the object domain. We illustrate this phenomenon in Fig.~\ref{fig:subsampledData}, where Fig.~\ref{fig:subsampledData}a corresponds to the Fourier transform of the original object, Fig.~\ref{fig:subsampledData}b corresponds to applying a direct inverse Fourier transform which recovers the original object, and Fig.~\ref{fig:subsampledData}c corresponds to applying an inverse Fourier transform to the cropped $64 \times 64$ part shown in Fig.~\ref{fig:subsampledData}a, which recovers a coarsened version of the object. 

For the object $\bfz$ and the gradient $\nabla_{\bfz} \Phi_\mathcal{M}$, we choose the restriction operator 
\begin{align}
  \bfI_h^H = \dfrac{1}{4}
    \left( 
    \begin{array}{@{}*{11}{c}@{}}
       1       & 1       & 0      & \hdots  &     0   & 1       & 1       & 0      &      0  &  \hdots  &    0     \\ 
       0       & 1       & 1      &     0   &  \hdots & 0       & 1       & 1      &    0   &  \hdots  &    0     \\
       \vdots  & \ddots  & \ddots & \ddots  &  \ddots &         & \ddots  & \ddots & \ddots &          & \vdots   \\
       0       & \hdots  & 0   &   1        &     1   & 0       &     0    & \hdots &    0   &     1    & 1 
    \end{array} 
    \right) \in \bbR^{n_H^2 \times n_h^2}
\end{align}
as described in \cite[pg. 69]{trottenberg2000multigrid}, which consists of a four-point average that maps from the vertices to the cell centers. The prolongation operator is chosen in the standard fashion as  
\begin{align}
  \bfI_H^h = c \: (\bfI_h^H)^\top,
\end{align}
for some constant $c>0$. As a standard choice for the cell-centered prolongation operator, we follow \cite[pg. 61]{trottenberg2000multigrid} and choose $c=4$ in our numerical experiments.
\subsection{Discussion}
MG/OPT is not an algorithm; instead, it is a framework that allows for different choices of discretization, restriction/prolongation operators, and optimization algorithms as its underlying solver. This flexibility is of particular importance in the context of ptychography as it allows us to exploit the natural hierarchy embedded in both the object domain and the data domain. Furthermore, since each MG/OPT cycle contains at least one iteration of the fine-grid optimization algorithm (OPT with $k_1+k_2>0$ in Sec.~\ref{subsec:MG/OPTcycle}), the linesearch in the coarse-grid correction guarantees the convergence of MG/OPT in the same fashion as OPT  \cite{nash2000multigrid}.
\begin{figure}[t]
  \centering
  \begin{tabular}{cc}
    \centering
    \includegraphics[width=0.45\textwidth, height=2.2in]{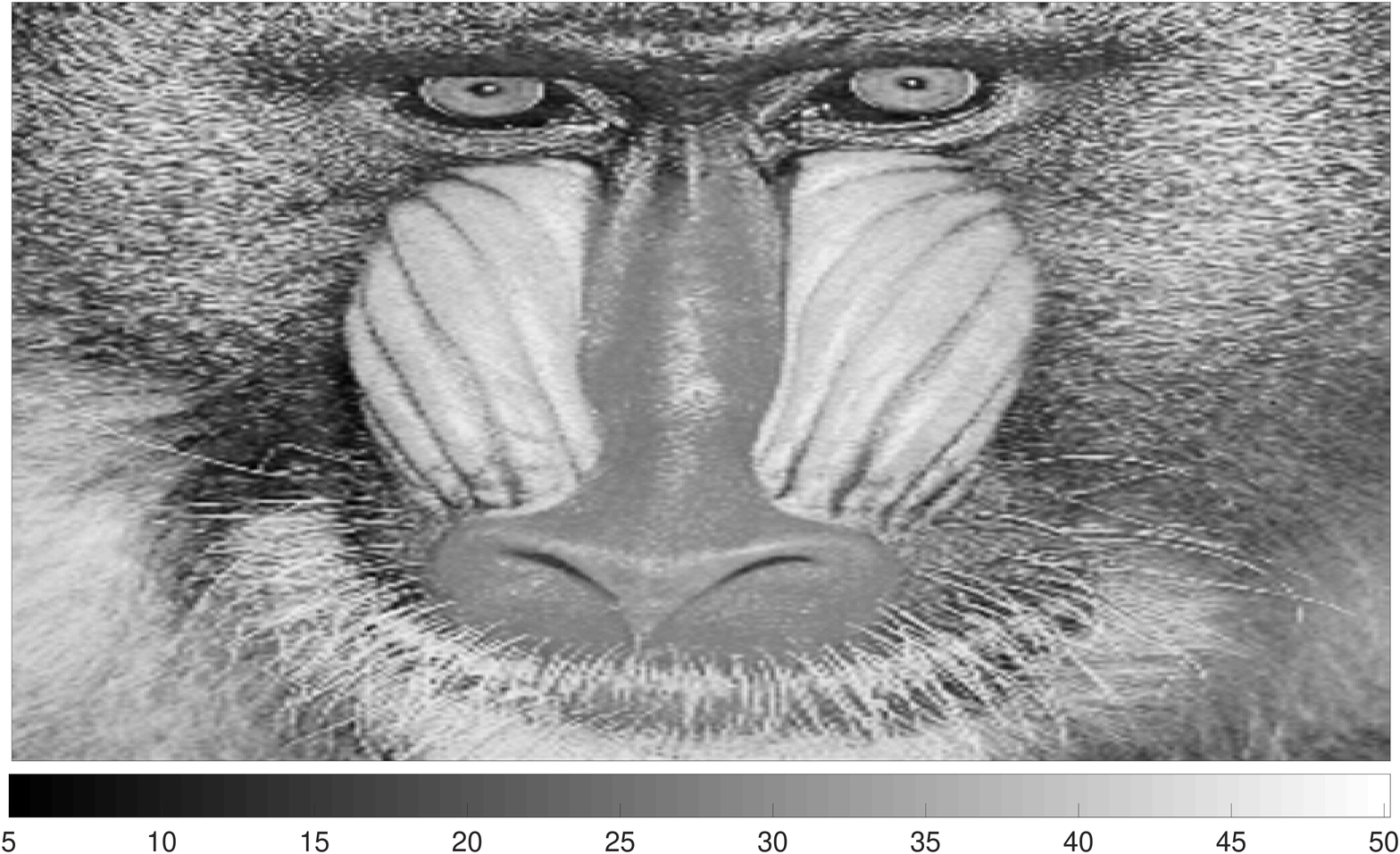}
    &
    \includegraphics[width=0.45\textwidth, height=2.2in]{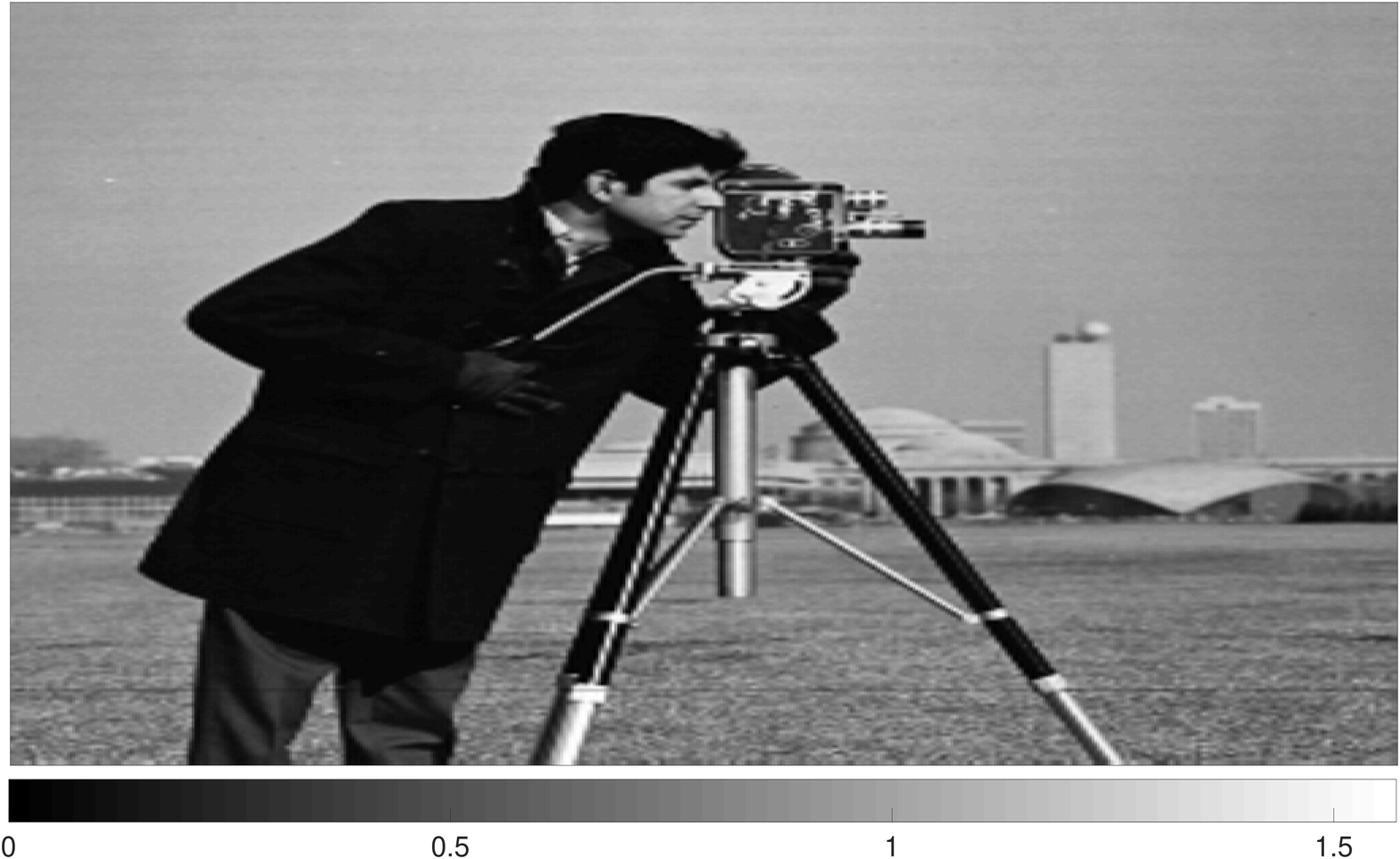}
  \end{tabular}
  \caption{\textit{Ground truth used to simulate data in numerical experiments. The baboon image is used as the magnitude and the camera man image is used as the phase of the object of interest.}}
  \vspace{-5mm}
  \label{fig:groundTruth}
\end{figure}
\begin{figure}[t]
  \setlength\tabcolsep{1 pt}
    \centering
    \begin{tabular}{ccccccccc}
        \includegraphics[width=0.105\textwidth, height=0.65in]{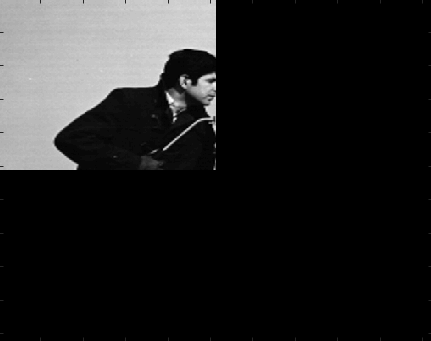}
        &
        \includegraphics[width=0.105\textwidth, height=0.65in]{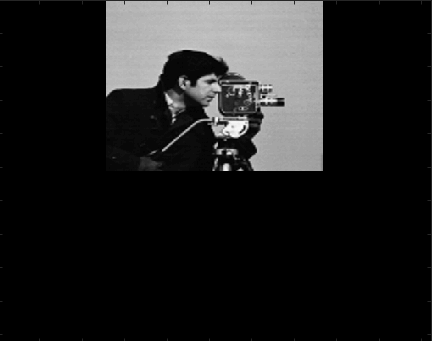}
        &
        \includegraphics[width=0.105\textwidth, height=0.65in]{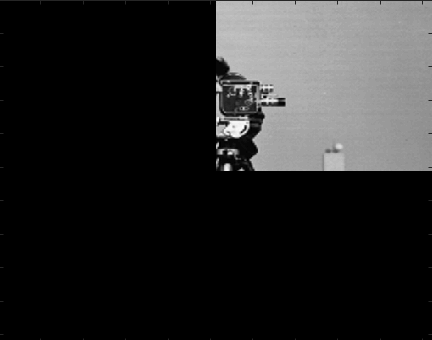}
        &
        \includegraphics[width=0.105\textwidth, height=0.65in]{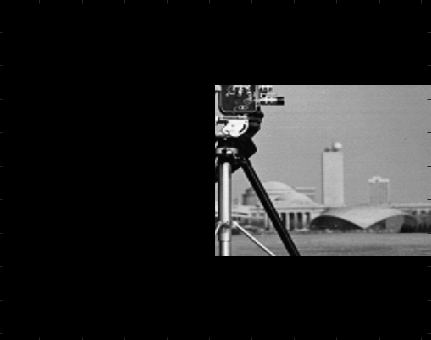}
        &
        \includegraphics[width=0.105\textwidth, height=0.65in]{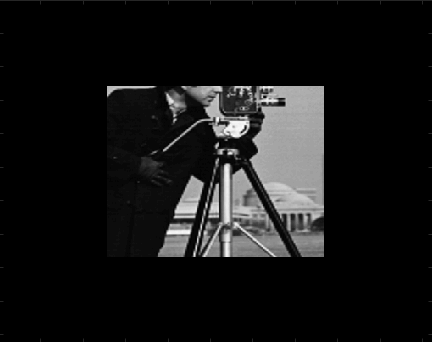}
        &
        \includegraphics[width=0.105\textwidth, height=0.65in]{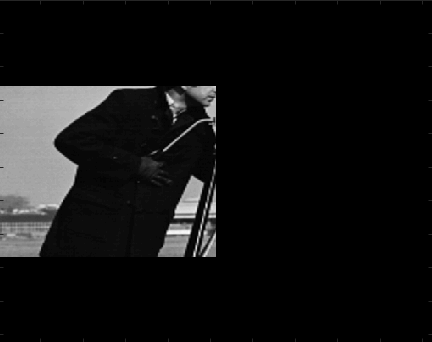}
        &
        \includegraphics[width=0.105\textwidth, height=0.65in]{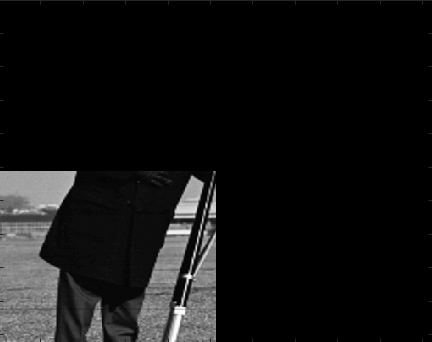}
        &
        \includegraphics[width=0.105\textwidth, height=0.65in]{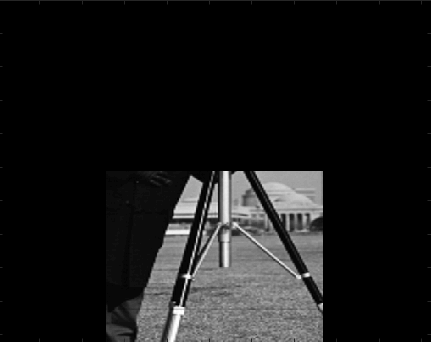}
        &
        \includegraphics[width=0.105\textwidth, height=0.65in]{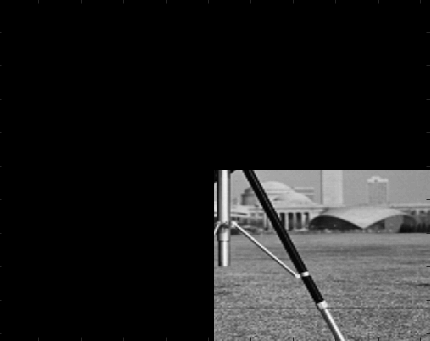}
        \\
        \small{$\bfQ_1$} & \small{$\bfQ_2$} & \small{$\bfQ_3$} & \small{$\bfQ_4$} & \small{$\bfQ_5$} & \small{$\bfQ_6$} & \small{$\bfQ_7$} & \small{$\bfQ_8$} & \small{$\bfQ_9$}
    \end{tabular}
    \caption{\textit{Scanning positions (probes) described in Sec.~\ref{subsec:setup}. Here, each probe illuminates a window of size $256 \times 256$ pixels and shifts $128$ pixels at a time, leading to $50\%$ overlap between consecutive probes.}}
    \label{fig:scanningProbes}
    \vspace{-5mm}
\end{figure}
\begin{figure}[t]
  \setlength\tabcolsep{1 pt}
  \centering
  \begin{tabular}{ccc}
    \includegraphics[width=0.3\textwidth, height=1.5in]{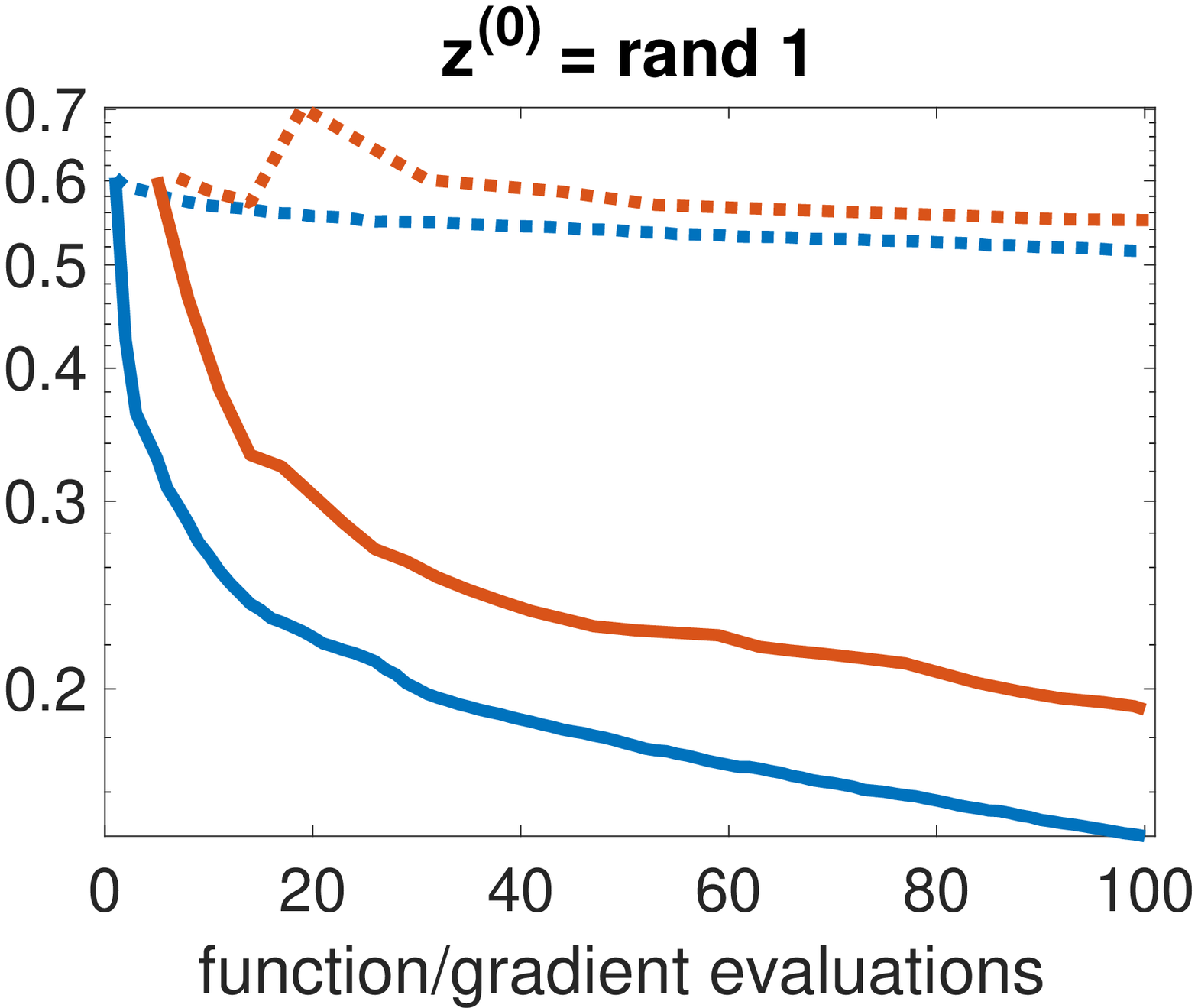}
    &
    \includegraphics[width=0.3\textwidth, height=1.5in]{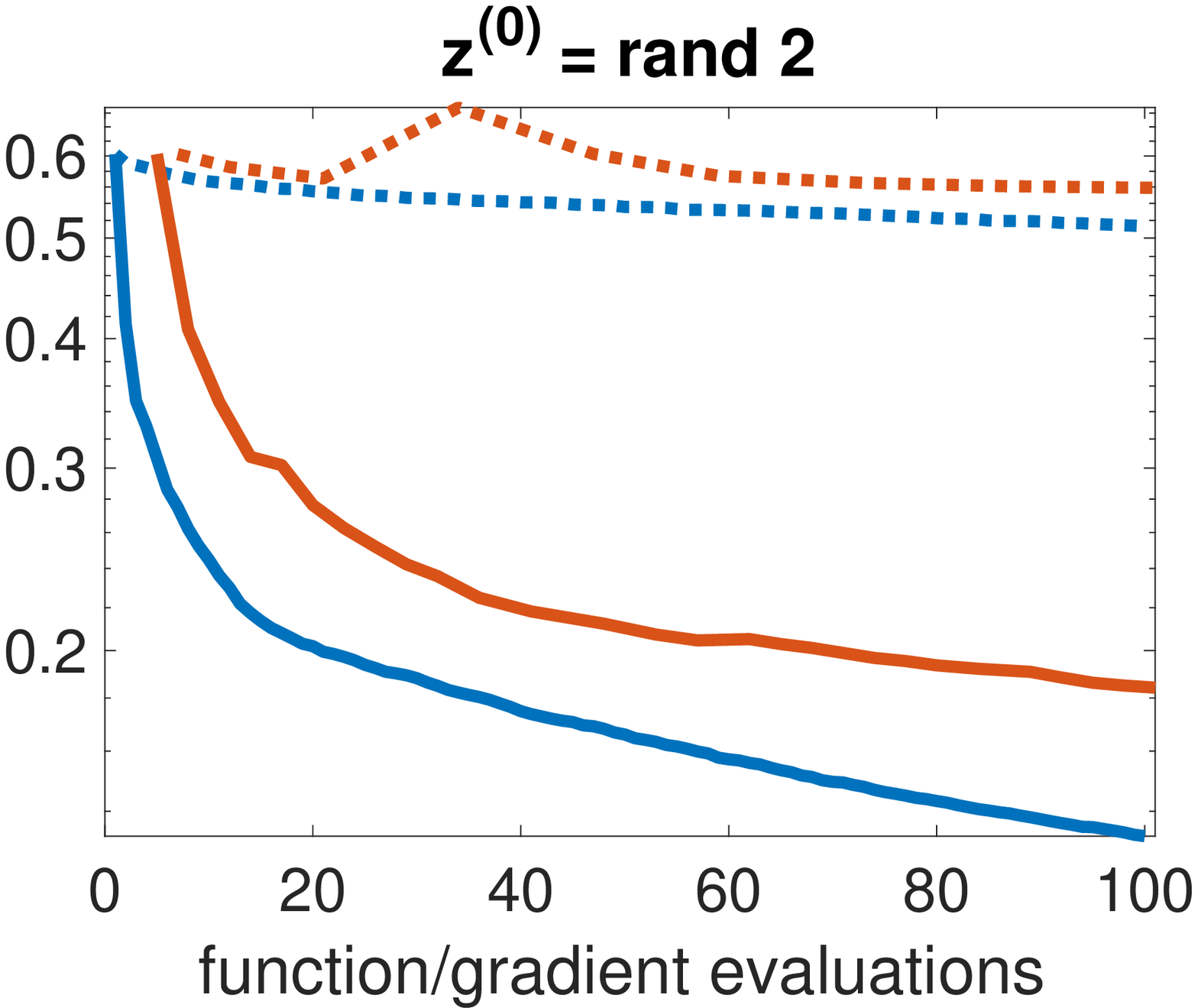}
    &
    \includegraphics[width=0.3\textwidth, height=1.5in]{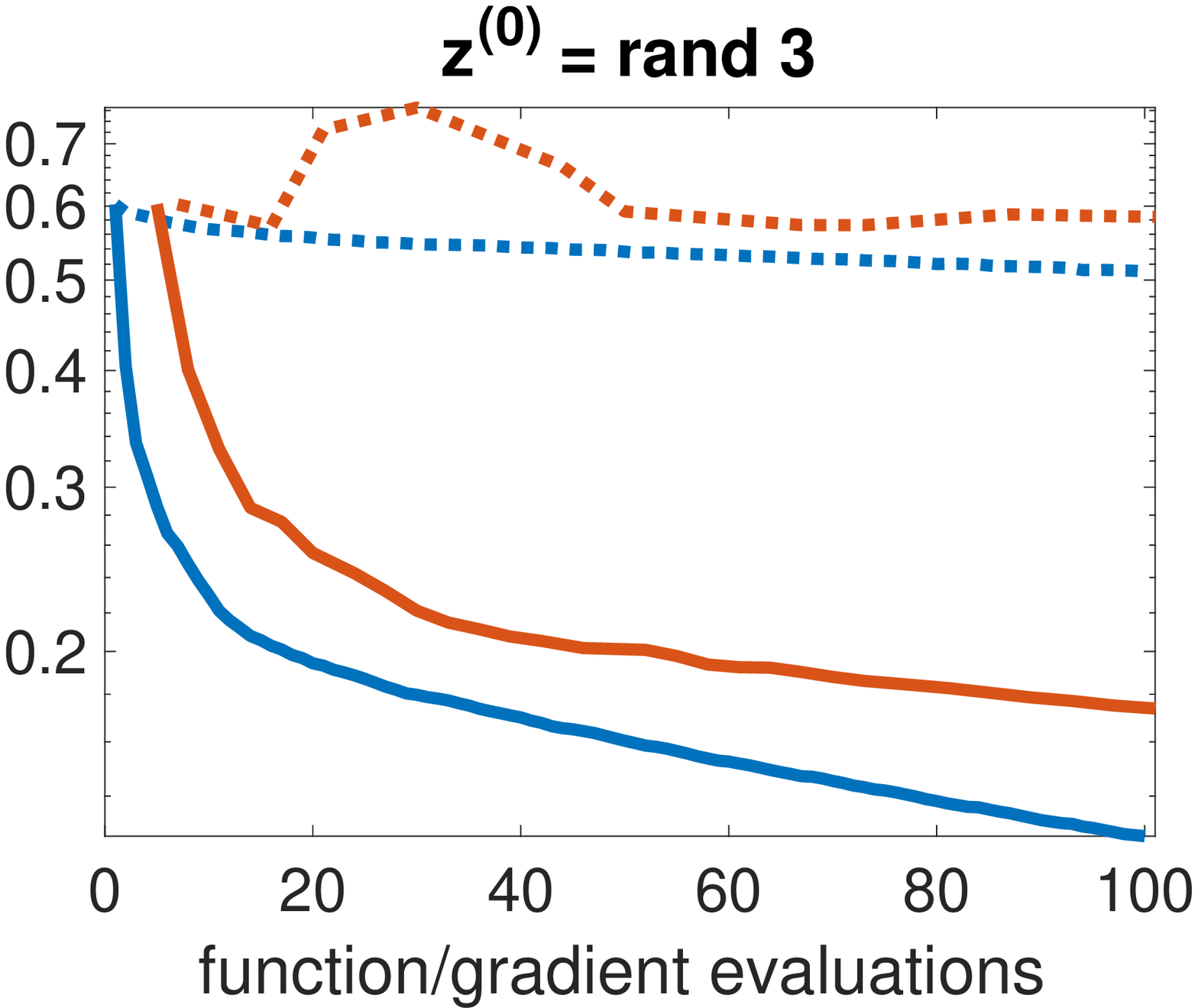}
    \\
    \includegraphics[width=0.3\textwidth, height=1.5in]{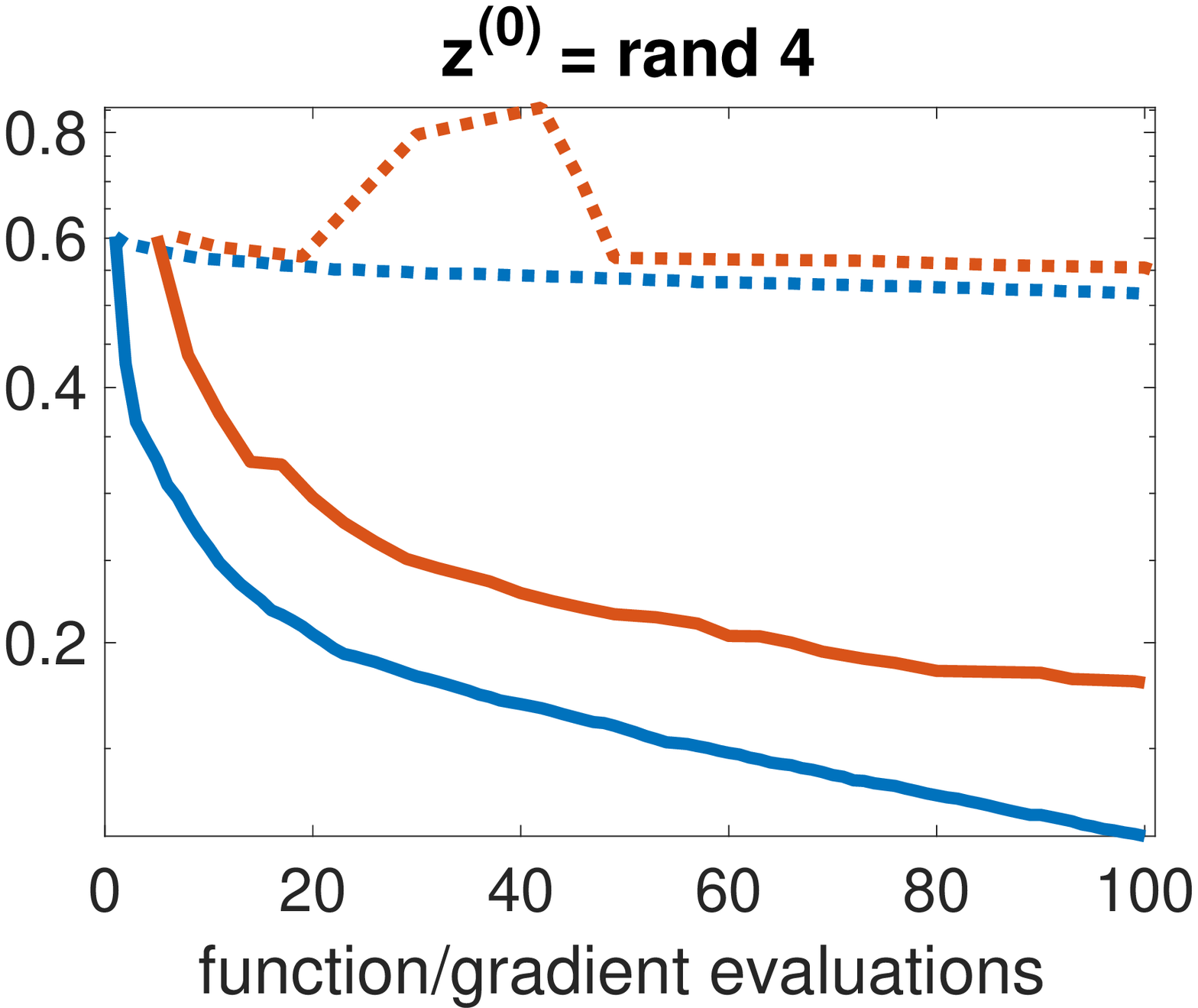}
    &
    \includegraphics[width=0.3\textwidth, height=1.5in]{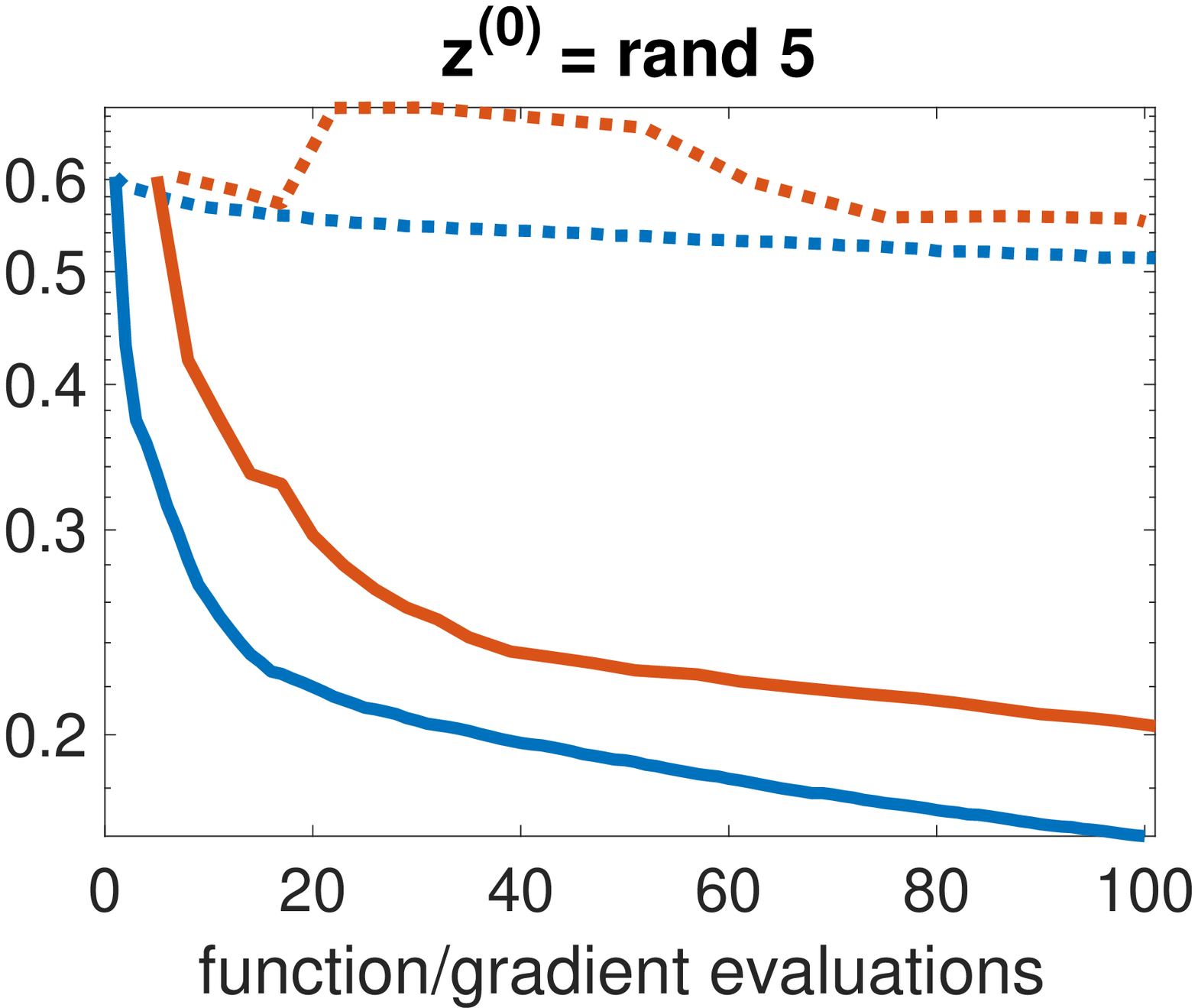}
    &
    \includegraphics[width=0.3\textwidth, height=1.5in]{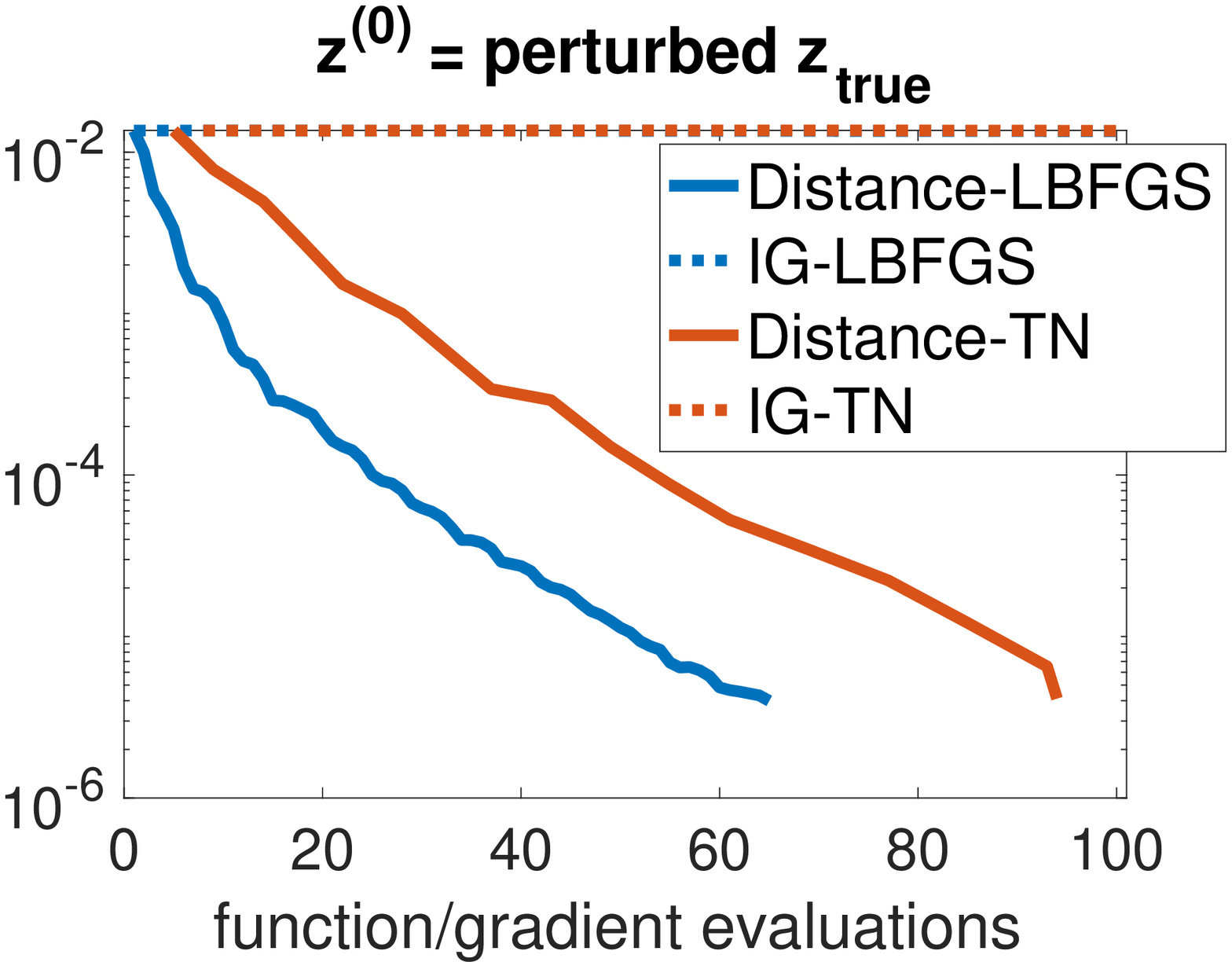}
\end{tabular}
  \caption{\textit{Relative errors for LBFGS and TN for two error metrics: the distance metric $\Phi_{\mathcal{M}}$ and Intensity Gaussian $\Phi_\mathcal{IG}$.}}
  \label{fig:compareMisfits}
  \vspace{-5mm}
\end{figure}
\section{Numerical Results}\label{sec:numResults}
In this section, we compare the performance of derivative-based algorithms for the two objective functions~\eqref{eq:misfit1} and~\eqref{eq:distanceMisfit} described in Sec~\ref{sec:mathFormulation}. We then demonstrate the potential of the MG/OPT scheme applied to a 2D ptychographic phase retrieval problem.
\subsection{Ptychography Setup}
\label{subsec:setup}
For our experiments, we consider a 2D ptychographic phase retrieval modeled by~\eqref{eq:forwardProblem}. We simulate our data using two $n \times n$ images, baboon and cameraman from MATLAB's demo images. We assume $m=n$, which is a practical assumption since the object resolution is typically determined by the data resolution. We set the baboon image to be the magnitude and the cameraman image to be the phase of our ground truth (see Fig.~\ref{fig:groundTruth}). In order to avoid the ambiguity of a global constant phase offset \cite{maiden2017further} (common in phase retrieval), we set the range of the true phase to be $[0,\frac{\pi}{2}]$. Consequently, we compute the phase errors in our experiments by first mapping the current iterates to this range before comparing with the true phase. For simplicity, we use binary probes, where each probe illuminates $\frac{n}{2} \times \frac{n}{2}$ pixels, that is, $\bfQ_k \in \mathbb{R}^{n^2 \times n^2}$ is a diagonal matrix that satisfies
\begin{align}
[\bfQ_k]_{j,j}= 
\begin{cases}
  1 & \text{if pixel}\quad j \quad \text{is illuminated} \\
  0 & \text{otherwise}.\\
\end{cases}
\end{align}
The illumination window is shifted $\frac{n}{4}$ pixels at a time starting from the top-left corner to the right until it reaches the top-right corner. The probe is then shifted downwards and then to the left until it reaches the left edge again. This scanning procedure is continued until we have covered the entire image, and results in a total of $9$ probes with $50\%$ overlap between adjacent probes. We illustrate this in Fig.~\ref{fig:scanningProbes} for all the different scanning positions. We note that the scanning patterns and overlap percentage in our setup are consistent with the recommended experimental setups for ptychography \cite{bunk2008influence,huang2014optimization}.

\subsection{Error Metric Comparison}
\label{subsec:compareMisfits}
We begin by comparing the performance of derivative-based optimization algorithms for the two error metrics: $\Phi_\mathcal{IG}$ and $\Phi_{\mathcal{M}}$. To this end, we employ the truncated Newton (TN) algorithm described in \cite{nash1985preconditioning}, which uses a preconditioned conjugate gradient to solve the Newton system and a finite difference scheme to approximate the action of the Hessian on a vector. We also employ an LBFGS algorithm that stores the 25 most recent vectors to approximate the inverse of the Hessian \cite{nocedal2006numerical}. For both solvers, we use a basic linesearch scheme that takes a step once a descent direction is computed, and set a maximum of $50$ linesearch iterations. We run the inversions for $6$ different initial guesses: 5 random guesses, and a good initial guess that consists of a perturbation of the true solution. 
In Fig.~\ref{fig:compareMisfits}, we show performance of LBFGS and TN applied to $\Phi_\mathcal{IG}$ and $\Phi_{\mathcal{M}}$, in terms the relative reconstruction error 
\begin{align} \label{eq:relErrFormula}
  \frac{\left\| \bfz^{(j)}-\bfz_{\rm true} \right\|}{\left\| \bfz_{\rm true} \right\|},
\end{align} where $\bfz_{\rm true}$ is the ground truth. For a fair comparison, the performance is reported based on the number of function/gradient evaluations.  



For all initial conditions, we observe the following. First, we obtain much higher accuracies when optimizing over $\Phi_{\mathcal{M}}$, since gradient-based algorithms applied to $\Phi_\mathcal{IG}$ tend to get stuck at potentially poor-quality local minima. 
This observation is consistent with the success of PIE and other alternating projection algorithms, since they optimize the more stable objective function $\Phi_{\mathcal{M}}$ for ptychography. This observation is also seen in other ptychography experiments \cite{qian2014efficient,yeh2015experimental}. Second, in terms of computational work, LBFGS outperforms TN for this particular application because each inner CG iteration in TN requires an additional gradient computation, which leads to an overall higher computational load; we discuss this in Sec~\ref{subsec:CompCosts}. Consequently, in the remainder of our experiments, we consider only the LBFGS as OPT to the distance metric $\Phi_\mathcal{M}$.
\subsection{Multilevel Ptychography}
We illustrate the potential of MG/OPT for the ptychographic phase retrieval problem described in Sec.~\ref{subsec:setup}.
\subsubsection{Computational Costs} \label{subsec:CompCosts}
We measure the computational cost of different algorithms by counting the number of fine-grid function/gradient evaluations. Each evaluation of $\Phi_{\mathcal{M}}$ requires $N$ projections $\mathcal{P}_{\mathcal{M}_k}$, and the residuals from $\Phi_\mathcal{M}$ can be used to compute the gradients at negligible costs (see~\eqref{eq:distanceGradient}). Moreover, each projection consists of one Fourier transform and one inverse Fourier transform, which have complexity $\mathcal{O}(n \log (n))$. Therefore, one function/gradient evaluation has a total complexity of $\mathcal{O}(Nn \log(n))$. Similarly, each outer PIE iteration also contains $N$ projections $\mathcal{P}_{\mathcal{M}_k}$ and therefore has the same complexity as one function/gradient evaluation of $\Phi_{\mathcal{M}}$.

To account for the computational cost of MG/OPT (with LBFGS as its underlying solver), we determine (for each level) the relative cost of a function/gradient evaluation compared to an evaluation on the finest grid. That is, given our restriction approach, the computational cost of one function/gradient evaluation on the coarser grid is roughly one fourth of the cost of one function/gradient evaluation on the previous finer grid. Following this relationship, the computational work shown in Table~\ref{tab:workAlloc} amounts to $100$ function/gradient evaluations for each MG/Opt scheme. 
\begin{table}[t]
\centering
  \begin{tabular}{|c|c|ccccc|}
    \hline
    \multicolumn{7}{|c|}{\textbf{function/gradient evaluations of MG/OPT across grids}}
    \\
    \hline
    & noise level& $n=512$ & $n=256$ & $n=128$ & $n=64$ & $n=32$
    \\
    \hline
    \multirow{3}{*}{2-grid MG/OPT}
    & 0\% & 68 & 137 &  &  &
    \\
    & 5\% & 68 & 130 &  &  & 
    \\
    & 10\% & 72 & 114 &  &  & 
    \\
    \hline
    \multirow{3}{*}{5-grid MG/OPT}
    & 0\% & 88 & 32 & 24 & 6 &409
    \\
    & 5\% & 82 & 55 & 46 & 6 & 337
    \\
    & 10\% & 78 & 66 & 62 & 9 & 486
    \\
    \hline
  \end{tabular}
\caption{\textit{Allocation of computational work for MG/OPT inversions across multiple grids for different noise levels. We show the different number of linesearch iterations performed on the different grids.}}
\label{tab:workAlloc}
\vspace{-5mm}
\end{table}
\subsubsection{Inversion Setup}
We solve the phase retrieval problem using four approaches: PIE, LBFGS, a 2-level, and 5-level MG/OPT with LBFGS as its underlying solver. In the latter three approaches, we use the same linesearch scheme described in Sec.~\ref{subsec:compareMisfits} with a maximum of $50$ iterations. In LBFGS, we store the $25$ most recent vectors used to approximate the inverse of the Hessian at each iteration. We run the inversions for $0\%$, $5\%$, and $10\%$ added Gaussian noise using a random initial guess, and stop after $100$ function/gradient evaluations have been performed. 

In the MG/OPT scheme, we set a maximum of $3$ and $100$ iterations at the coarsest grids of the 2-level and 5-level schemes per cycle, respectively. We do this since an iteration on the coarsest grid of a 5-level V-cycle is much cheaper than an iteration on the coarsest grid of a 2-level V-cycle (see Sec~\ref{subsec:CompCosts}). As for the remaining finer grids, we set a maximum of $1$ iteration, i.e., $k_1=k_2=1$. For the inner stopping criteria, we stop when the relative gradient norm is less than $10^{-4}$ at the coarsest grid, and $10^{-3}$ at the remaining finer grids.
In the 2-level MG/OPT, we project once from $\bbC^{512 \times 512}$ onto $\bbC^{256 \times 256}$. In the 5-level MG/OPT, we project four times: $\bbC^{512 \times 512} \mapsto \bbC^{256 \times 256} \mapsto \bbC^{128 \times 128} \mapsto \bbC^{64 \times 64} \mapsto \bbC^{32 \times 32}$.
\newcommand{\rottexttt}[1]{\rotatebox{90}{\hbox to 42mm{\hss #1\hss}}}
\begin{figure}[t]
    \setlength\tabcolsep{0 pt}
    \centering
    \begin{tabular}{cccc}
      & \textbf{\small{$\mathbf{0\%}$ noise}} & \textbf{\small{$\mathbf{5\%}$ noise}}& \textbf{\small{$\mathbf{10\%}$ noise}}
      \vspace{-5mm}
      \\
      \rottexttt{\small{$\Phi_{\mathcal{M}}$}}
      &
      \includegraphics[width=0.3\textwidth, height=1.5in]{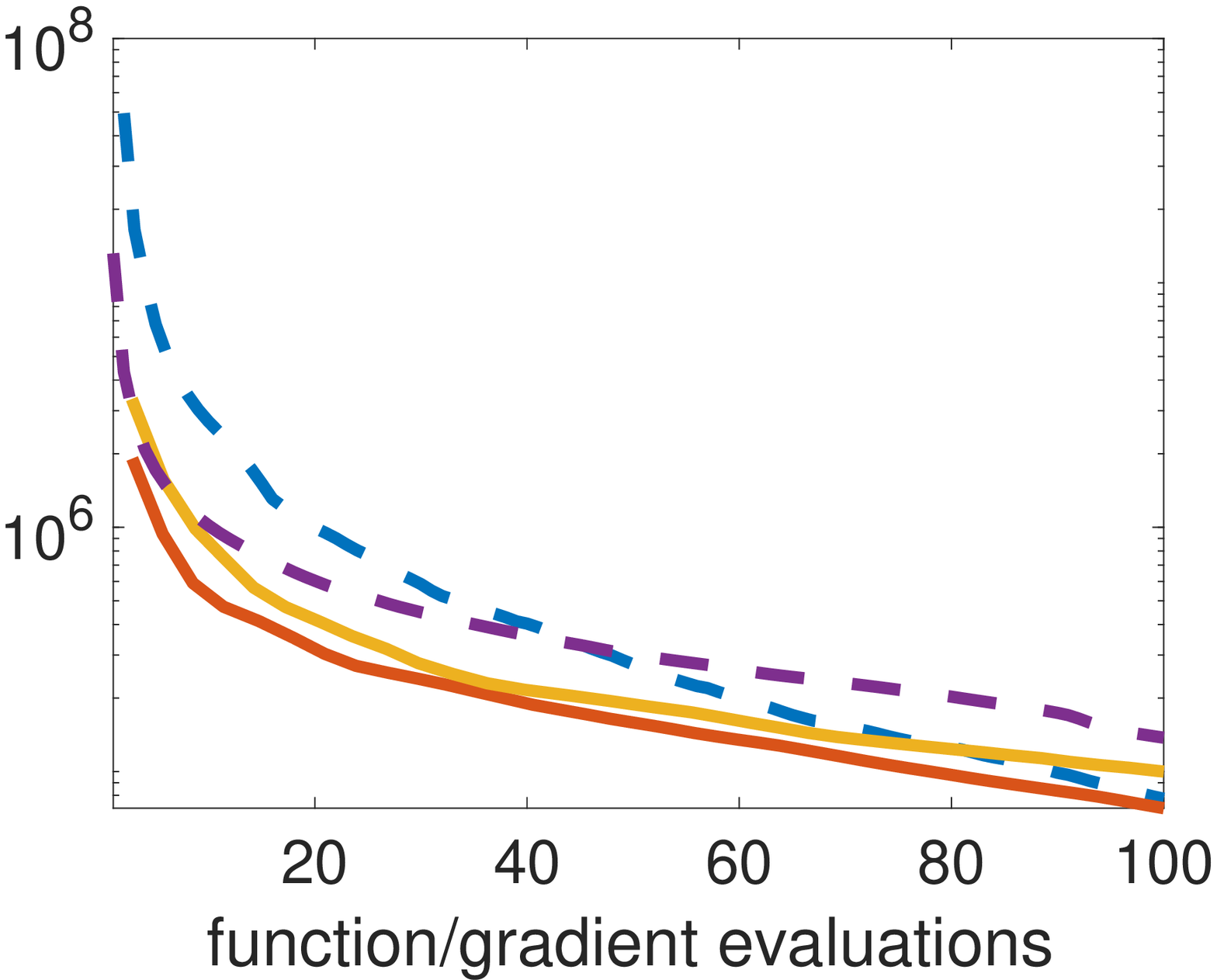}
      &
      \includegraphics[width=0.3\textwidth, height=1.5in]{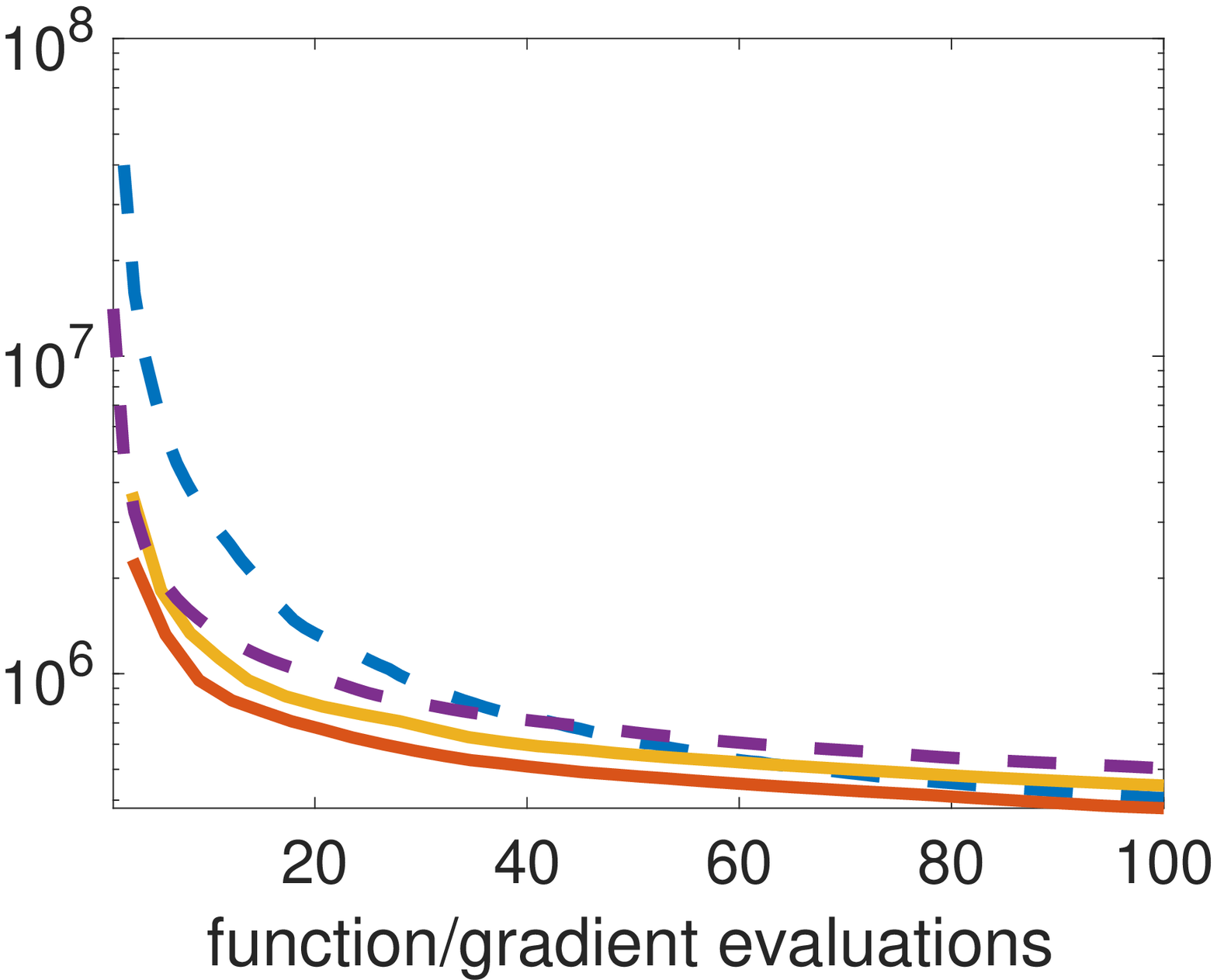}
      &
      \includegraphics[width=0.3\textwidth, height=1.5in]{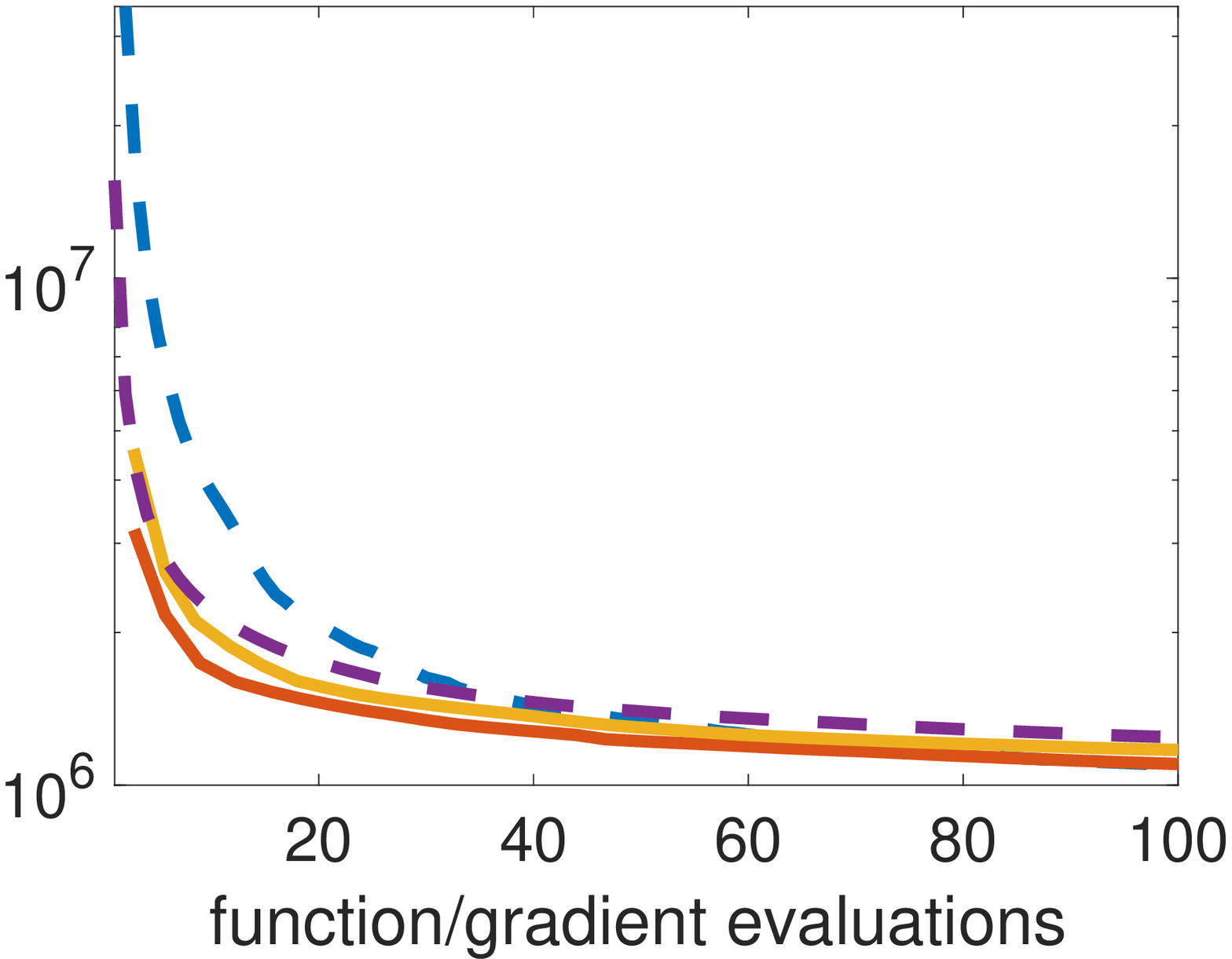}
      \\
      \rottexttt{\small{relative errors}}
      &
      \includegraphics[width=0.3\textwidth, height=1.5in]{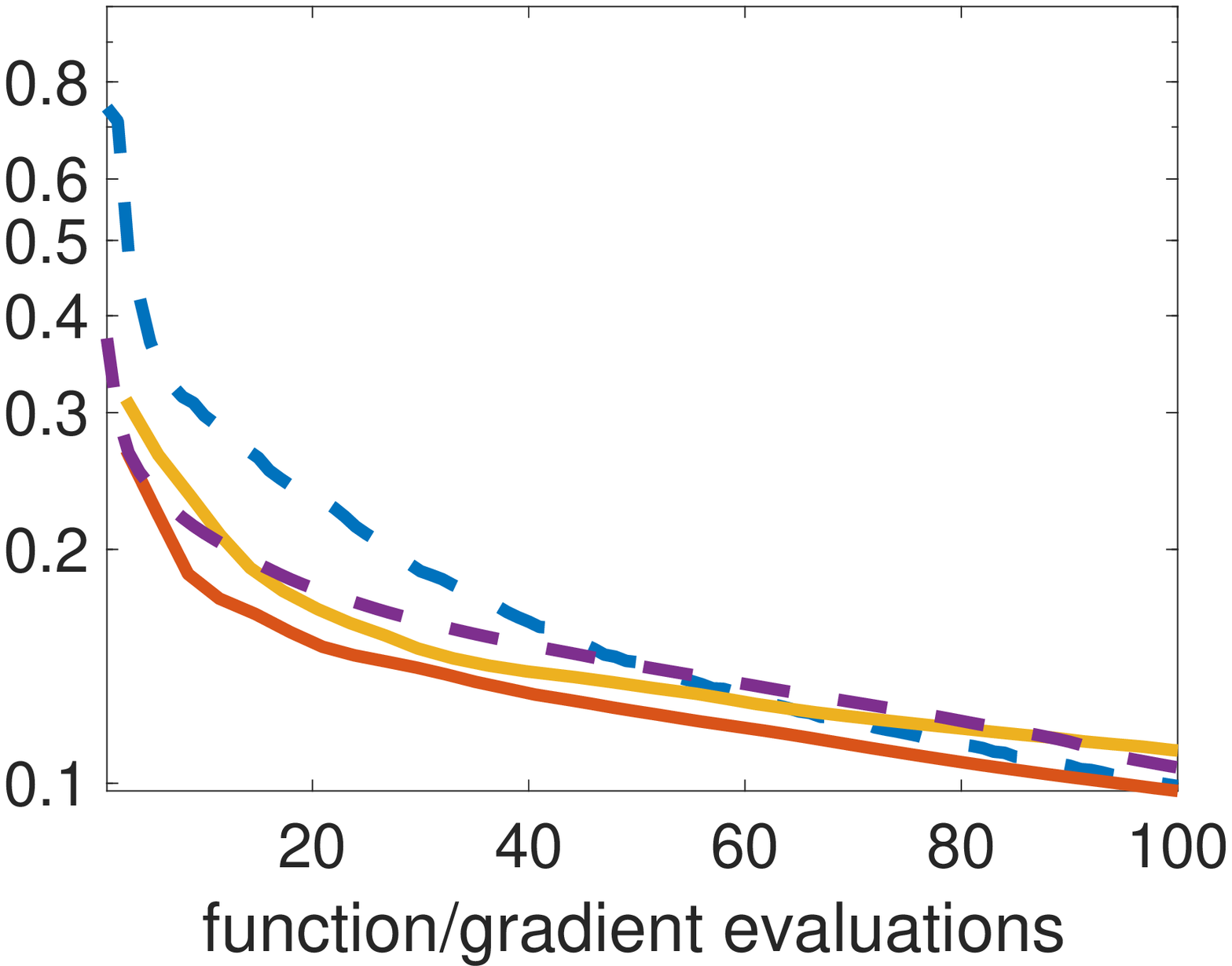}
      &
      \includegraphics[width=0.3\textwidth, height=1.5in]{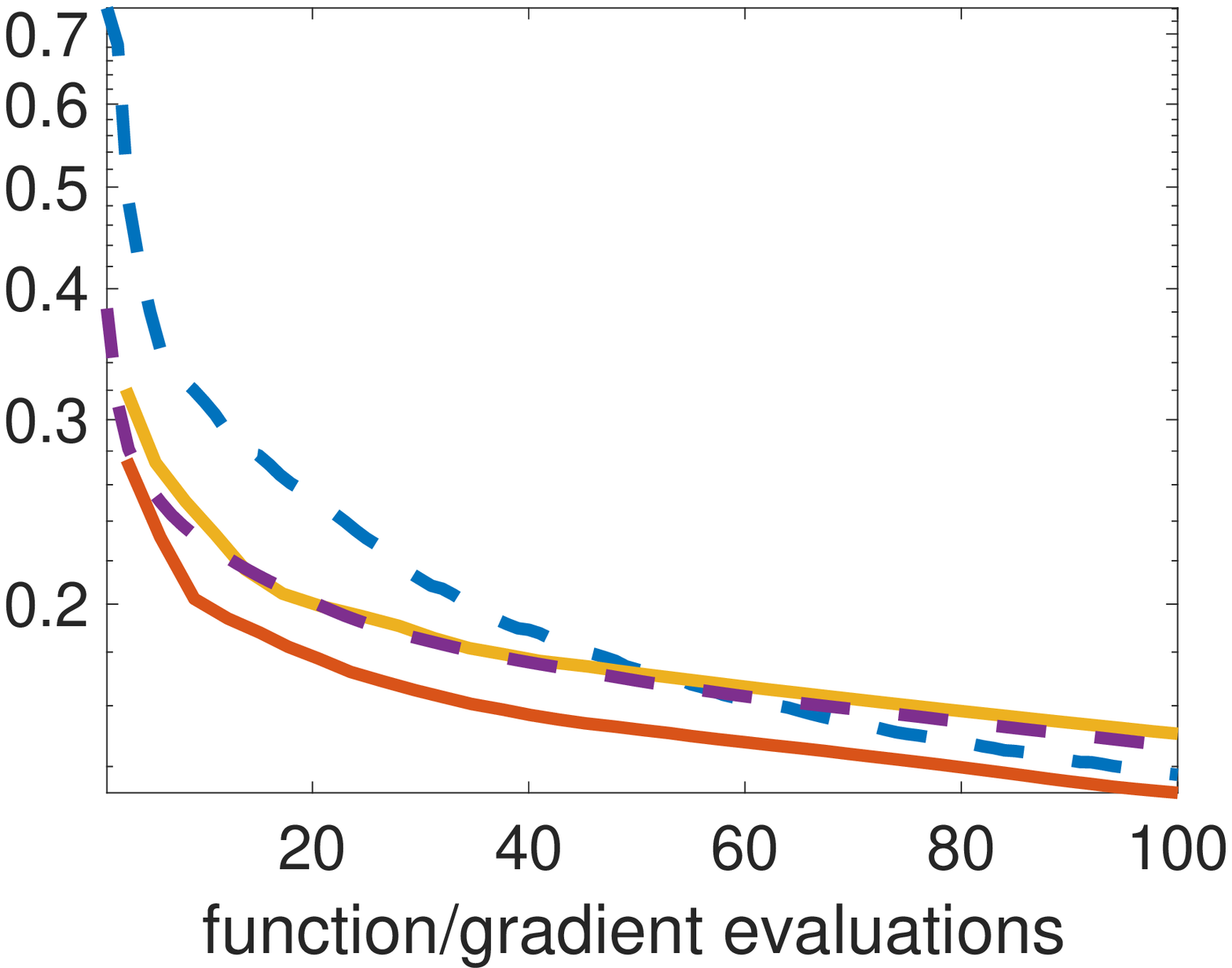}
      &
      \includegraphics[width=0.3\textwidth, height=1.5in]{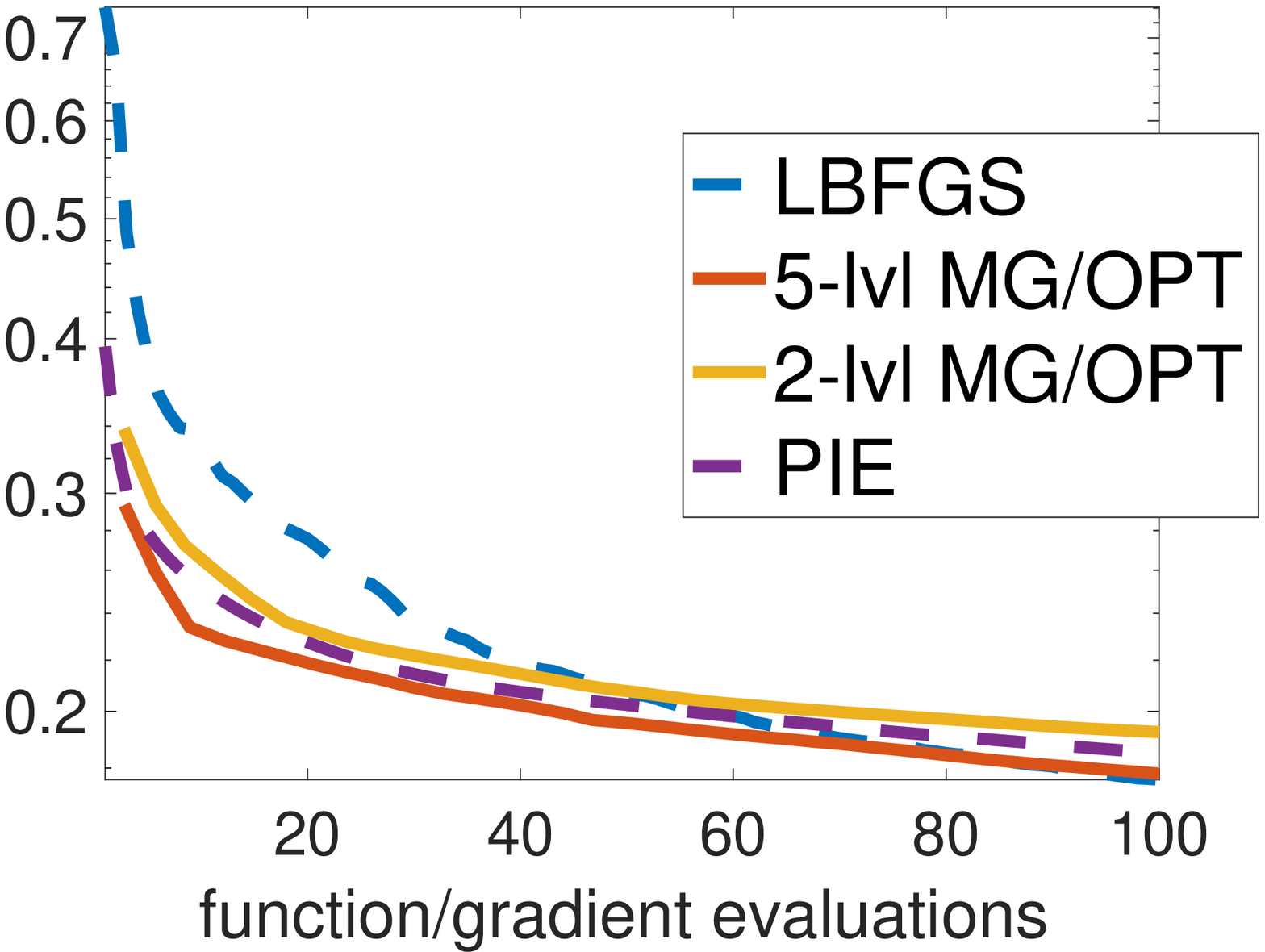}
    \end{tabular}
    \caption{\textit{Objective function and relative error histories for 0\%, 5\%, and 10\% noise levels.}}
    \label{fig:convHis}
    \vspace{-5mm}
\end{figure}

\begin{figure}[t]
    \setlength\tabcolsep{0 pt}
    \centering
    \begin{tabular}{cccc}
      & \textbf{\small{$\mathbf{0\%}$ noise}} & \textbf{\small{$\mathbf{5\%}$ noise}}& \textbf{\small{$\mathbf{10\%}$ noise}}
      \vspace{-5mm}
      \\
      \rottexttt{\small{magnitude errors}}
      &
      \includegraphics[width=0.3\textwidth, height=1.5in]{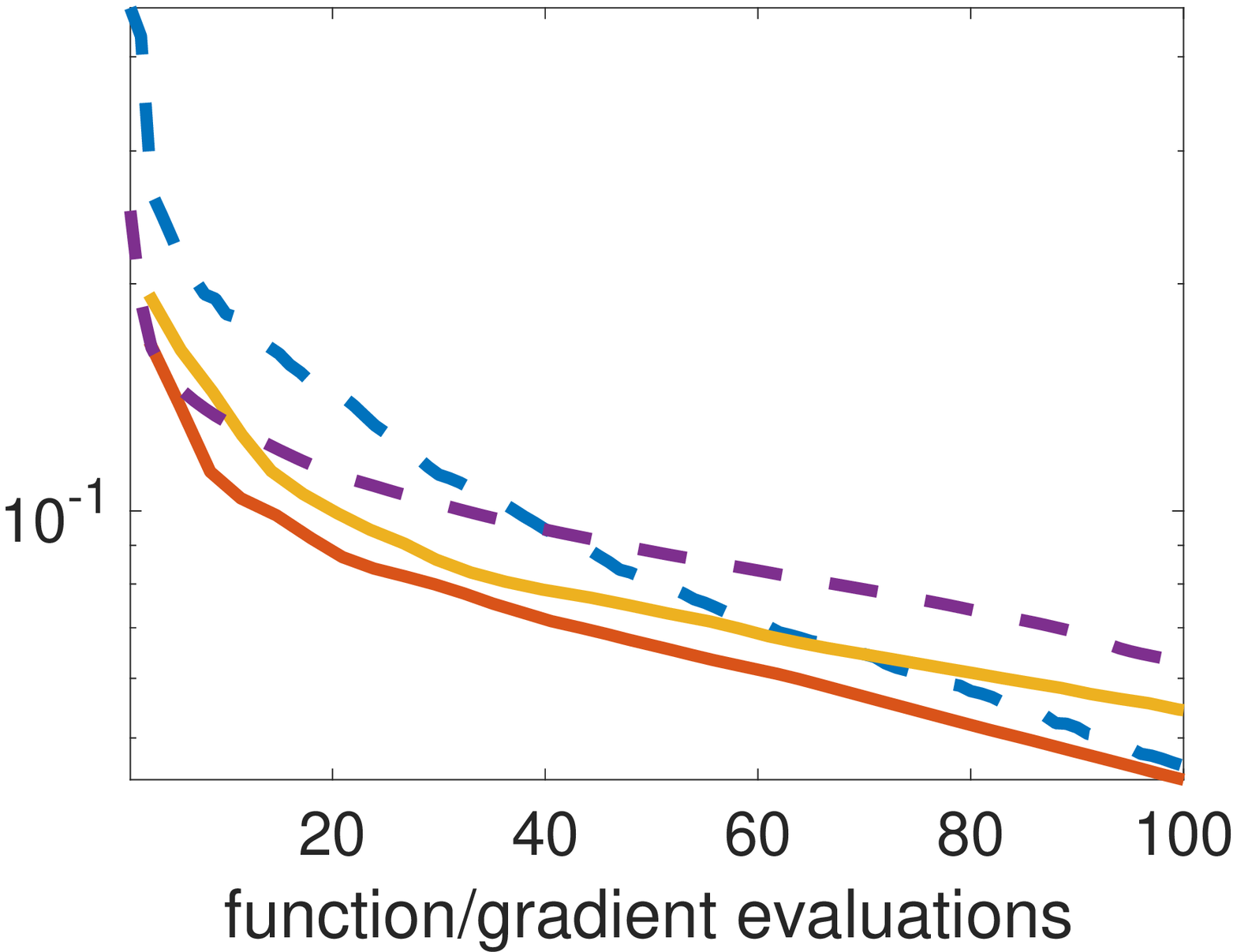}
      &
      \includegraphics[width=0.3\textwidth, height=1.5in]{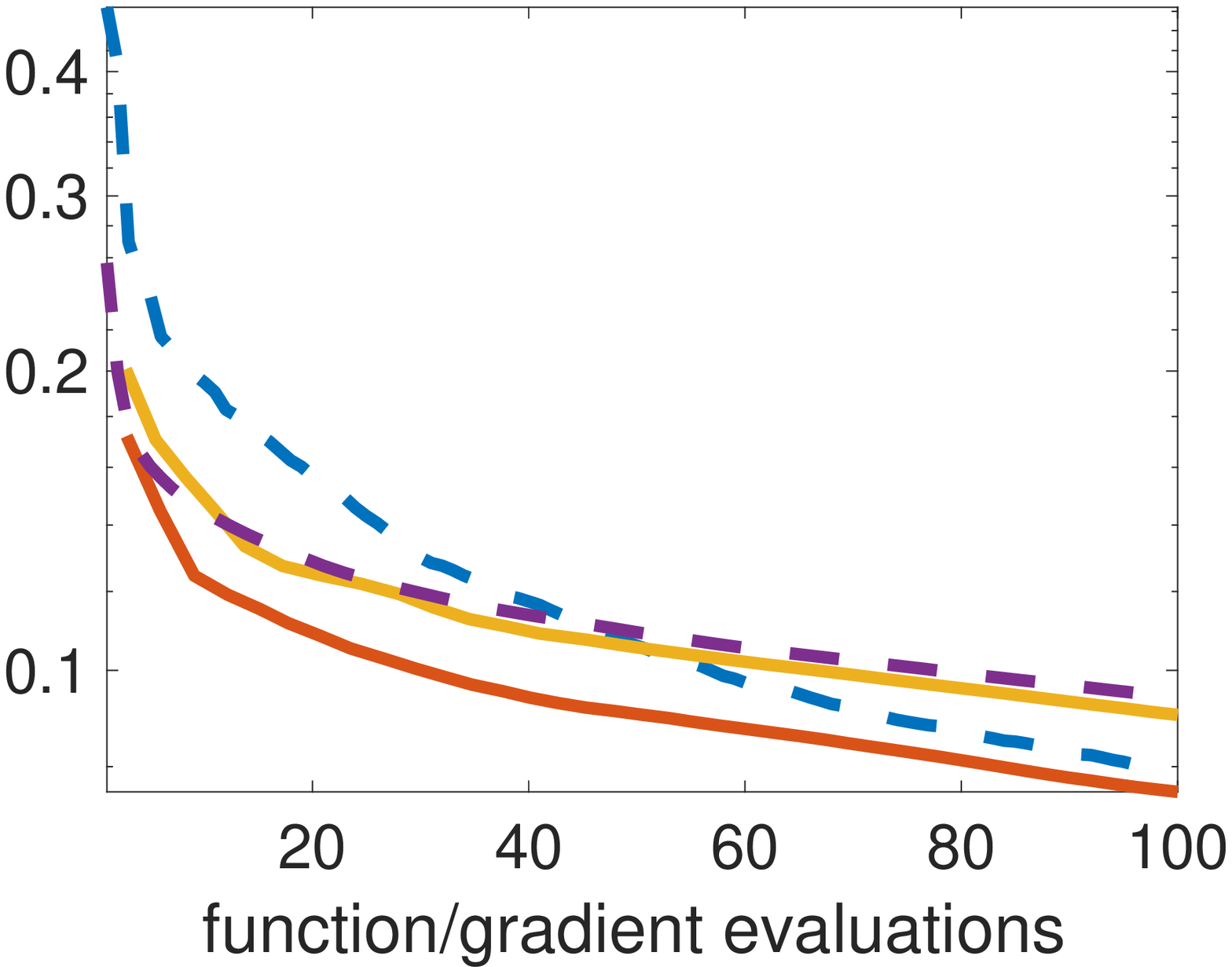}
      &
      \includegraphics[width=0.3\textwidth, height=1.5in]{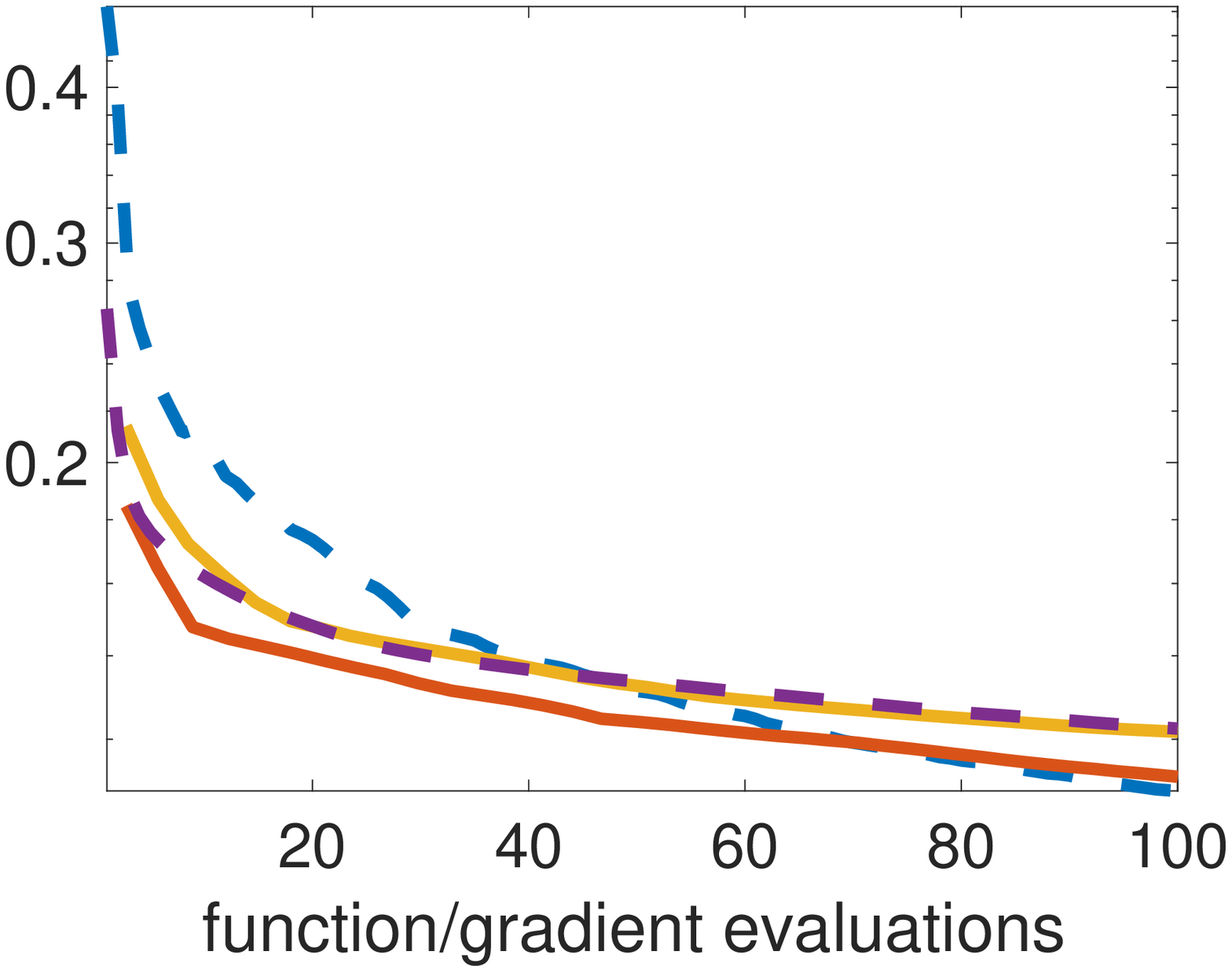}
      \\
      \rottexttt{\small{phase SSIM}}
      &
      \includegraphics[width=0.3\textwidth, height=1.5in]{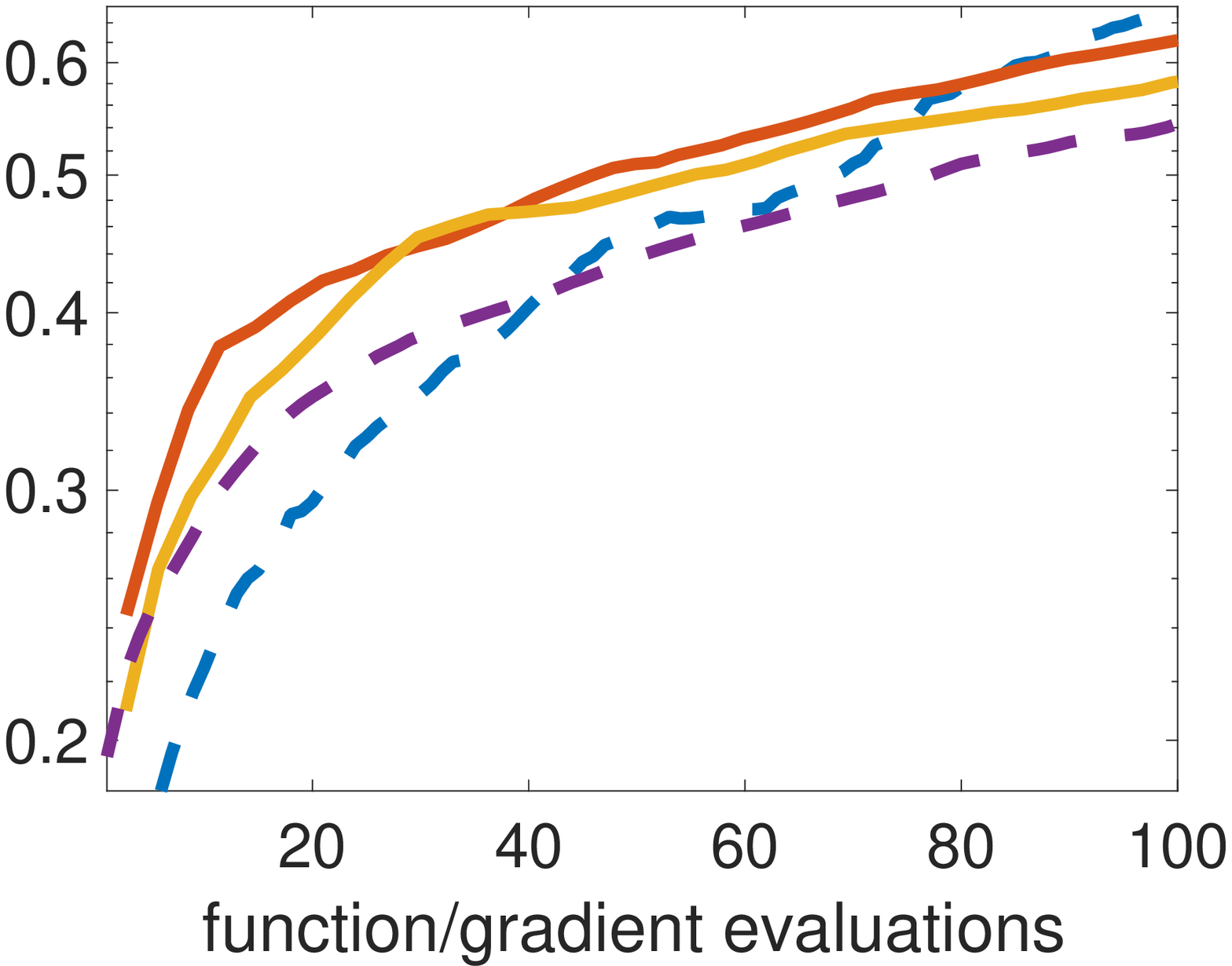}
      &
      \includegraphics[width=0.3\textwidth, height=1.5in]{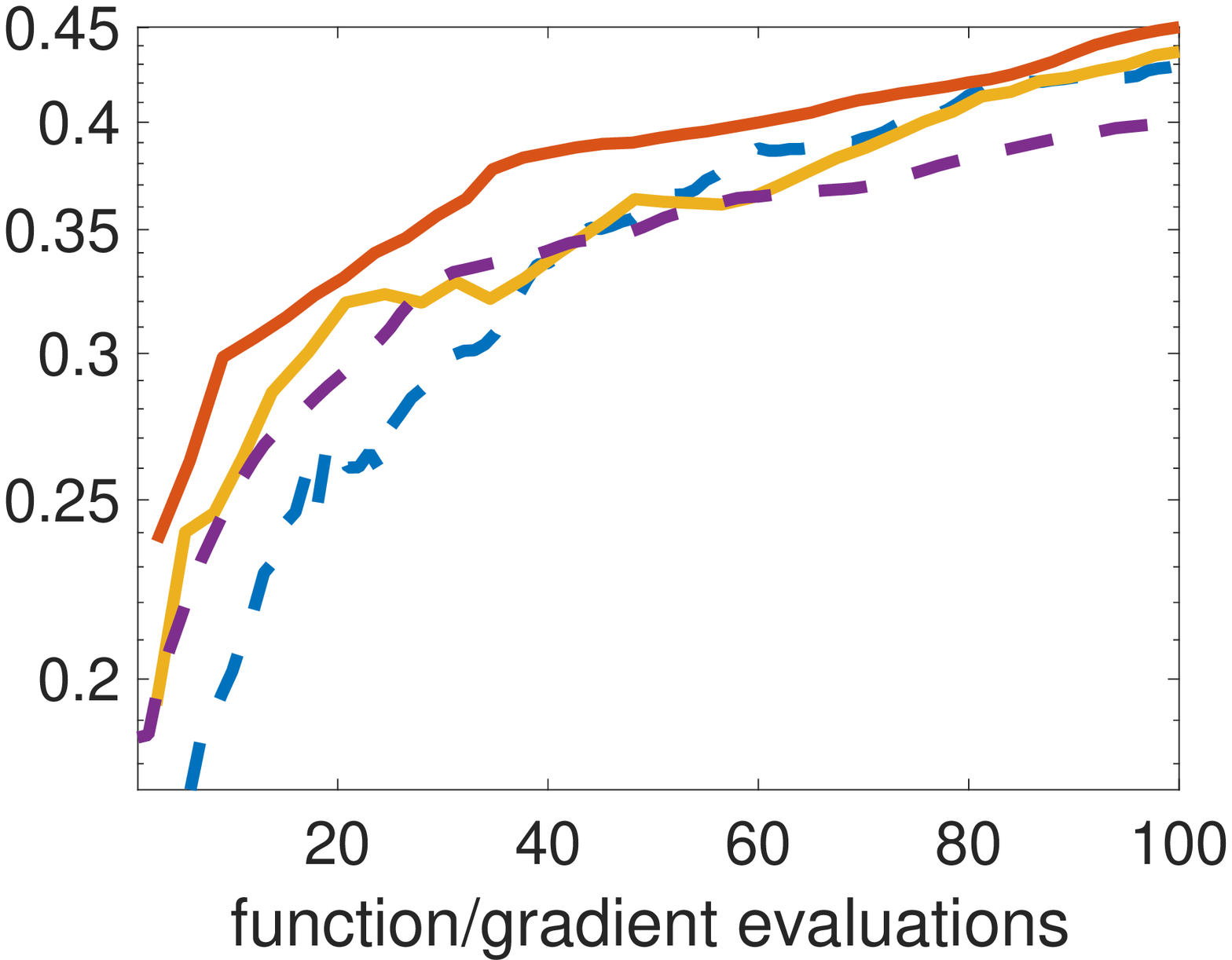}
      &
      \includegraphics[width=0.3\textwidth, height=1.5in]{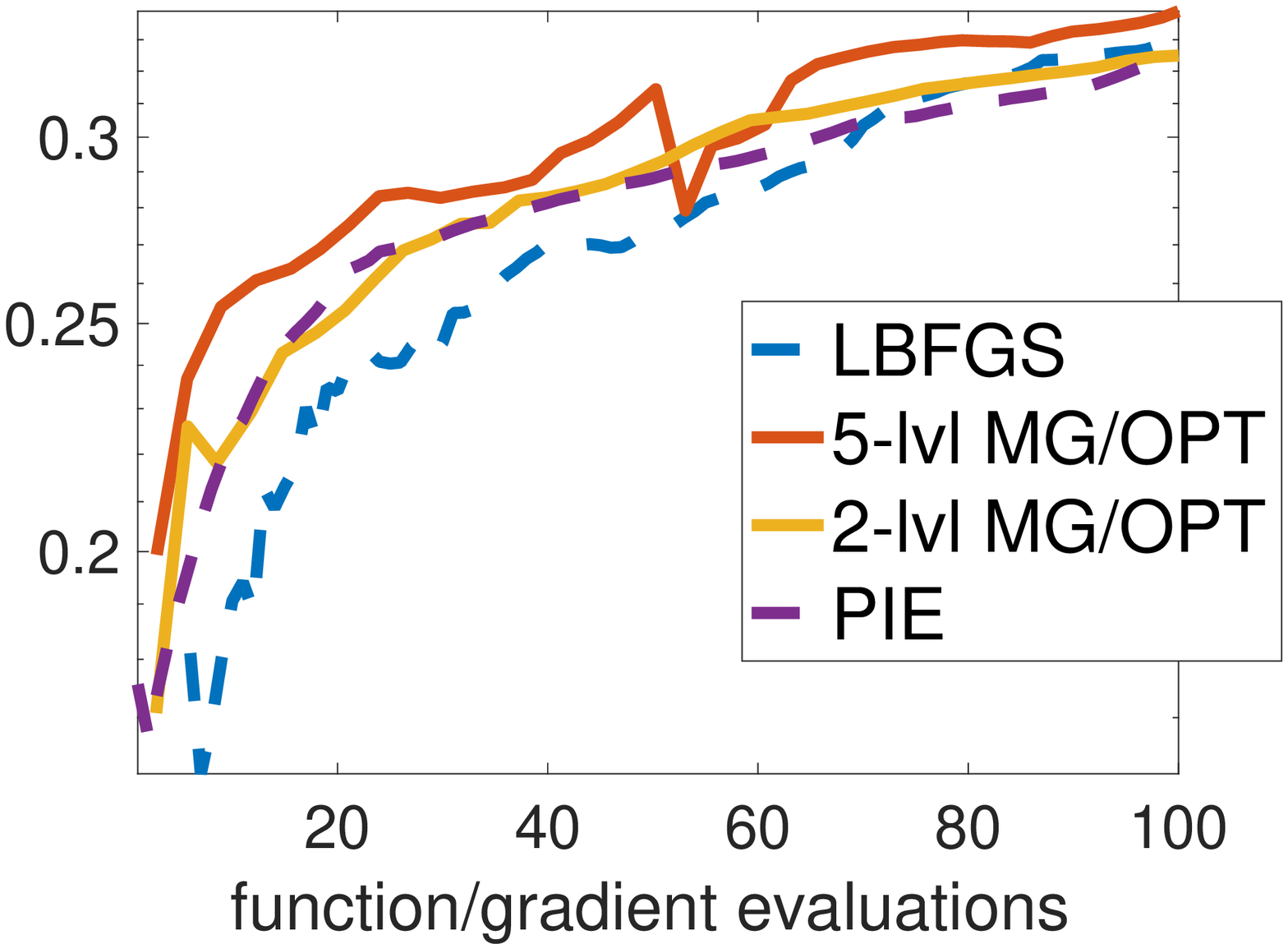}
    \end{tabular}
    \caption{\textit{Relative magnitude error and phase SSIM histories for 0\%, 5\%, and 10\% noise levels.}}
    \label{fig:phaseHis}
    \vspace{-5mm}
\end{figure}
\begin{table}[t]
\centering
\begin{tabular}{|c|c|c|c|c|}
\hline
& \multicolumn{2}{c|}{\textbf{convergence factor}} & \multicolumn{2}{c|}{\textbf{number of cycles}} 
\\
\hline
n & LBFGS & MG/OPT & LBFGS & MG/OPT 
\\
\hline
$64$ & 0.49 & 0.26 & 35 & 9 
\\
\hline
$128$ & 0.74 & 0.40 & 33 & 11 
\\
\hline
$256$ & 0.87 & 0.63 & 64 & 19 
\\
\hline
$512$ & 0.90 & 0.69 & 86 & 22 
\\
\hline
\end{tabular}
\caption{\textit{Comparison of convergence factors between LBFGS and MG/OPT. We run LBFGS and MG/OPT until convergence. In MG/OPT, the grid is always projected onto a mesh of size $32 \times 32$.}}
\label{tab:convFactors}
\vspace{-5mm}
\end{table}
\subsubsection{Results}
In Fig.~\ref{fig:convHis}, we show the convergence plots of the objective function values and relative errors of the three algorithms for $0\%$, $5\%$, and $10\%$ noise levels. As expected, MG/OPT accelerates the convergence of LBFGS as a result of a substantial allocation of work to the coarser grids (see Table~\ref{tab:workAlloc} for work allocation). In our experiments, we observe that the coarse grids contribute mainly in the earlier iterations, when there are coarse feature errors. After the first few cycles, the stopping criteria in the coarser grids are satisfied, and most of the work is then shifted to the finer grids. 
This behavior also leads to a lack of mesh-independence, where the convergence factor defined as
\begin{align} 
    c = \left(\dfrac{\Phi_{\mathcal{M}}(\bfz^{(j)}) - \Phi_{\mathcal{M}}(\bfz_{true})}{\Phi_{\mathcal{M}}(\bfz^{(0)}) - \Phi_{\mathcal{M}}(\bfz_{true})}\right)^{1/(j+1)}
\end{align}
is not fixed as the size of the problem grows. This can be seen in Table~\ref{tab:convFactors}, where we show these estimates as the size of the problem increases for a stopping criteria of the relative gradient norm $({\|\nabla \Phi_\mathcal{M} (\bfz^{j})\|}/{\|\nabla \Phi_\mathcal{M} (\bfz^{0})\|}) \leq 10^{-3}$.
Nonetheless, the contribution from the coarser grids causes the convergence factors of MG/OPT to be less sensitive to the problem size than the single-level LBGS. 

To complement Fig.~\ref{fig:convHis}, we show the reconstruction history of the magnitude and phase, separately, in Fig.~\ref{fig:phaseHis}. We use the same relative error metric shown in~\eqref{eq:relErrFormula} for the magnitudes. For the phase, however, we employ the structural similarity index metric (SSIM) \cite{wang2004image} as its error metric to further emphasize the perception change.
In particular, we measure the phase error by first performing a linear map of the current phase iterate to $[0,\frac{\pi}{2} ]$ (the domain of the true phase), and then perform the SSIM with the true phase as the reference image. An SSIM value of $1$ means the images are identical, and an SSIM value of 0 means that there are no structural similarities. The SSIM is an appropriate choice since it measures the image degradation as perceived change in structural information; this accounts for any potential global phase offset in our reconstructions.

As expected, the 5-level MG/OPT scheme outperforms the 2-level MG/OPT schemes, as we obtain higher computational savings coming from a more distributed allocation of work across grids (see Table~\ref{tab:workAlloc}). 
Fig.~\ref{fig:convHis} and Fig.~\ref{fig:phaseHis} also show that MG/OPT with LBFGS as its underlying solver outperforms PIE, a workhorse in the optics community.

In Fig.~\ref{fig:noiseRobustness}, we demonstrate the robustness of MG/OPT by plotting the relative errors vs. the noise levels after 38, 78, and 100 function/gradient evaluations. As can be seen in the plot, the 2-level and 5-level MG/OPT schemes are more robust to noise. We note that we are only able to record the values of the relative errors in the MG/OPT schemes after every V-cycle. Moreover, since the number of function/gradient evaluations performed after each V-cycle in MG/OPT is arbitrary (considering the relationship of the computational costs across grids described in Sec.~\ref{subsec:CompCosts}), we are unable to plot the relative errors after exactly 38, 78, or 100 function/gradient evaluations. Instead, we find an approximate function/gradient iteration value where all the MG/OPT schemes best agree (in our case, after 38, 78, and 100 function/gradient evaluations). For instance, the relative error shown in Fig.~\ref{fig:noiseRobustness}b for the 5-level MG/OPT at $10\%$ noise is recorded after 28 V-cycles, which in this particular case is equivalent to having performed $77.4$ function/gradient evaluations.

In Fig.~\ref{fig:reconstructions}, we show the ptychographic phase retrieval reconstructions for PIE, LBFGS, the 2-level MG/OPT, and the 5-level MG/OPT schemes for $0\%$, $5\%$ and $10\%$ noise levels. To allow for a fair comparison, the color axis is chosen identically to the ones shown in Fig.~\ref{fig:groundTruth}. The robustness of MG/OPT with respect to noise can be seen clearly in the phase reconstructions.
\newcommand{\rottext}[1]{\rotatebox{90}{\hbox to 27mm{\hss #1\hss}}}
\newcommand{\rottextt}[1]{\rotatebox{90}{\hbox to 4mm{\hss #1\hss}}}
\newcommand{\rottextttt}[1]{\rotatebox{90}{\hbox to 40mm{\hss #1\hss}}}
\begin{figure}[t]
    \setlength\tabcolsep{1 pt}
    \centering
    \begin{tabular}{cccc}
      & \small{a) 38 func/grad evals} & \small{b) 78 func/grad evals} & \small{c) 100 func/grad evals}
      \vspace{-2mm}
      \\
      \rottextttt{relative errors}
      &
      \includegraphics[width=0.3\textwidth, height=1.5in]{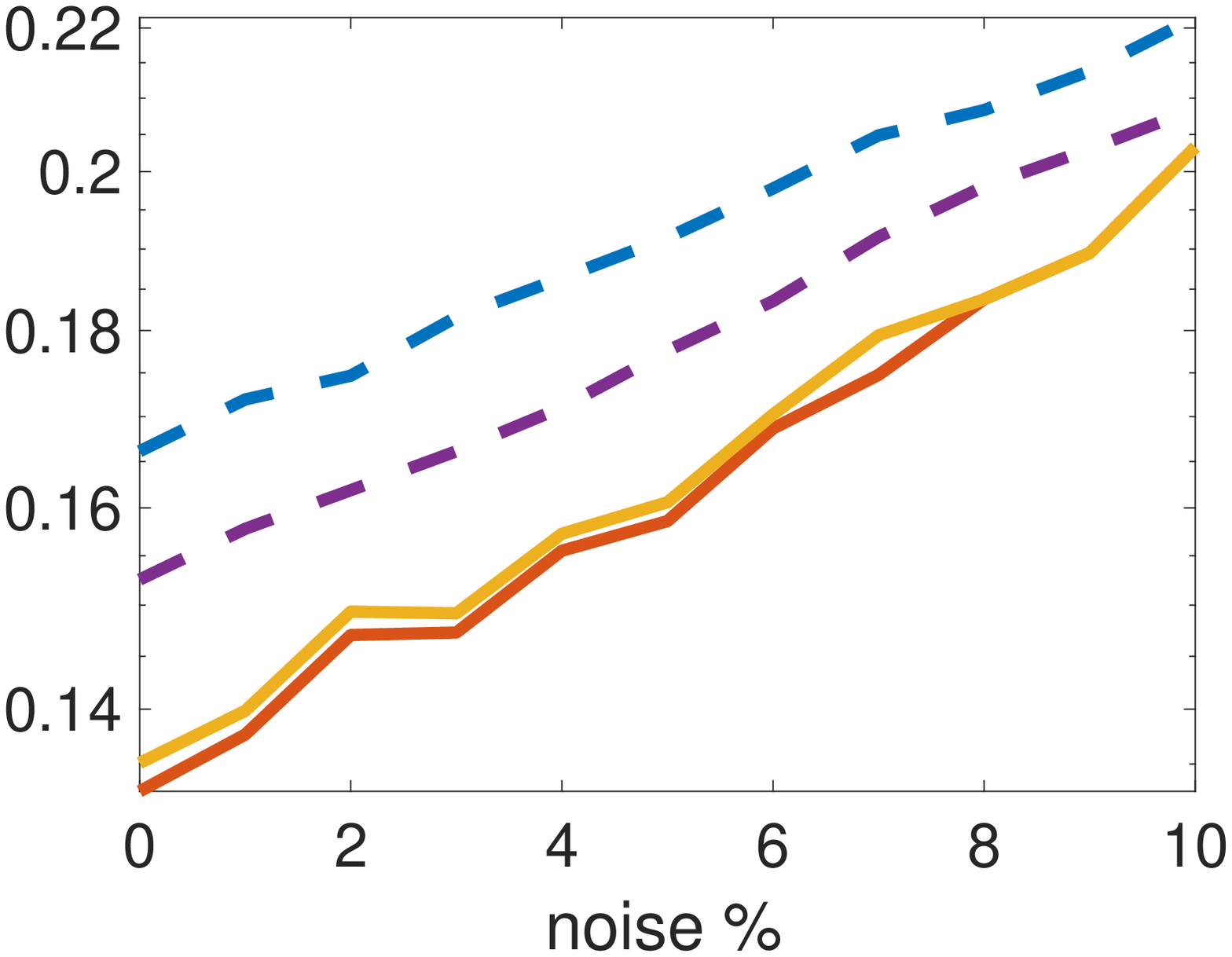}
      &
      \includegraphics[width=0.3\textwidth, height=1.5in]{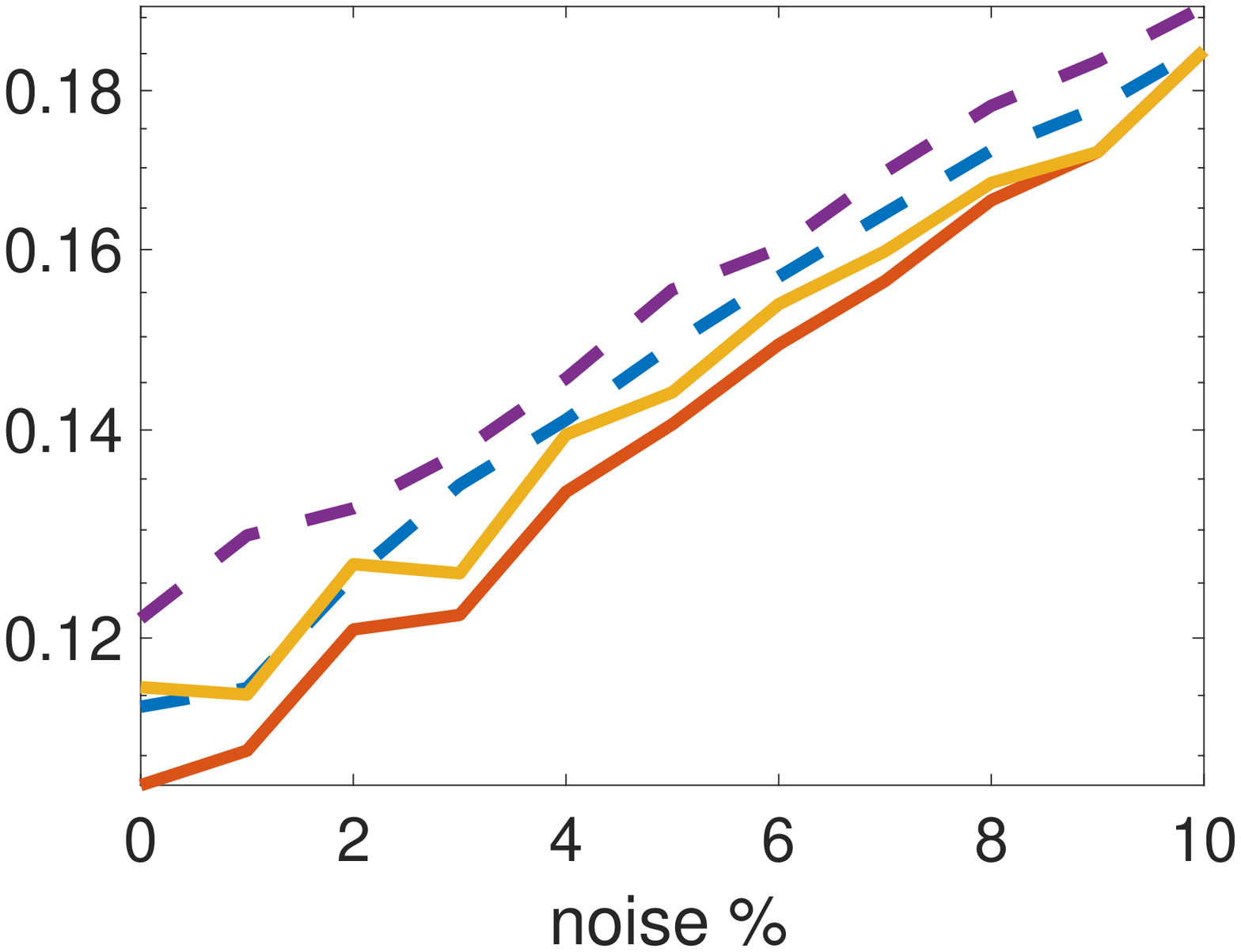}
      &
      \includegraphics[width=0.3\textwidth, height=1.5in]{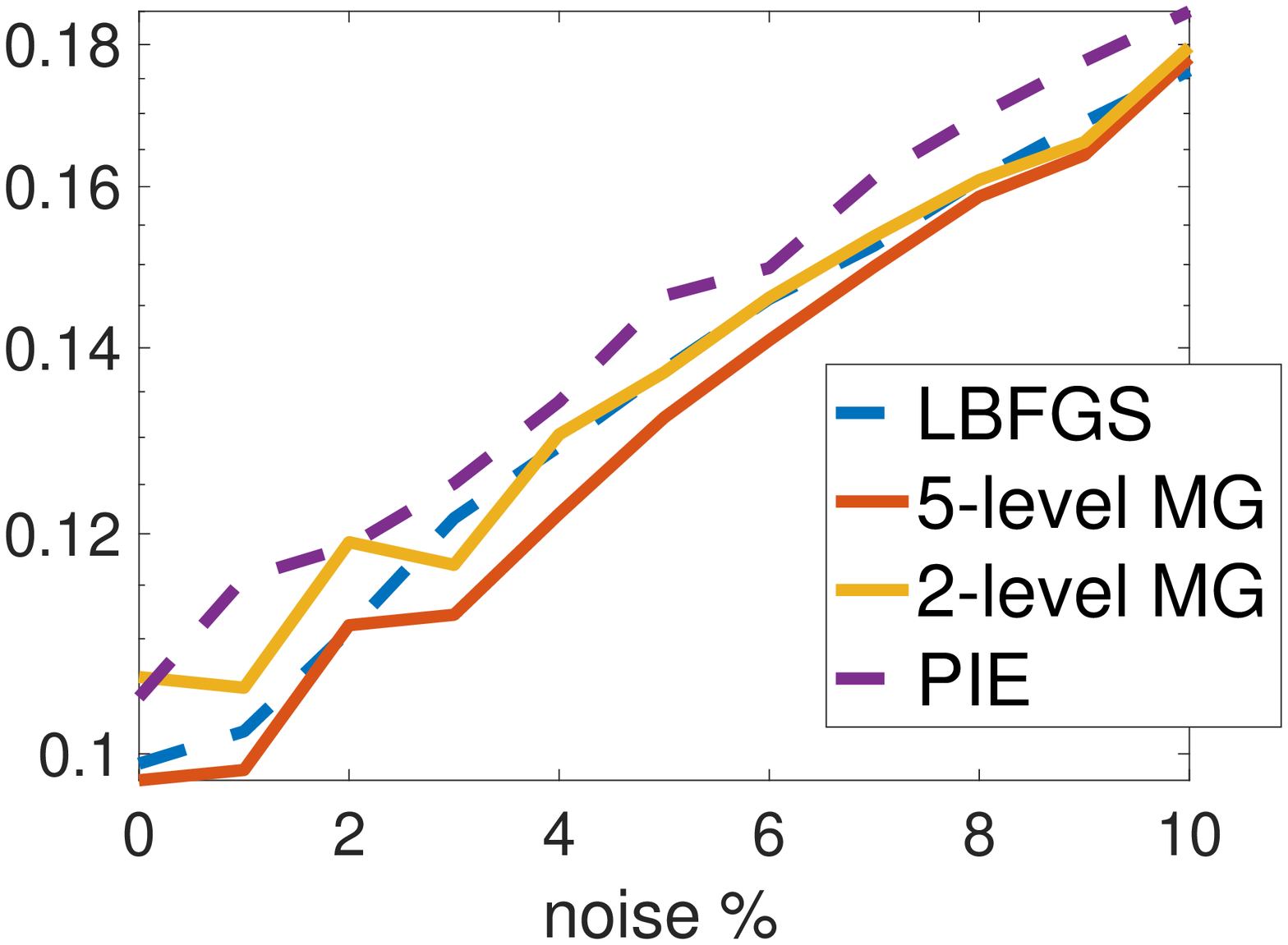}
    \end{tabular}
    \caption{\textit{Relative errors vs. noise after 38, 78, and 100 function/gradient evaluations.}}
    \label{fig:noiseRobustness}
\end{figure}
\begin{figure}[t]
 \setlength\tabcolsep{4 pt}
    \centering
    \begin{tabular}{cccccc}
    & &\textbf{LBFGS} & \textbf{PIE} & \textbf{2-lvl MG/OPT} & \textbf{5-lvl MG/OPT}
    \\
    \hline
    \multirow{2}{*}{\rottextt{\textbf{$\mathbf{0\%}$ noise}}}
    &
    \rottext{magnitude}
    &
    \includegraphics[width=0.21\textwidth]{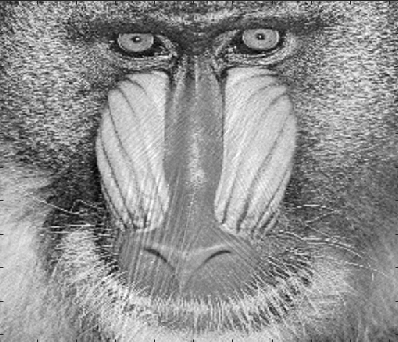}
    &
    \includegraphics[width=0.21\textwidth]{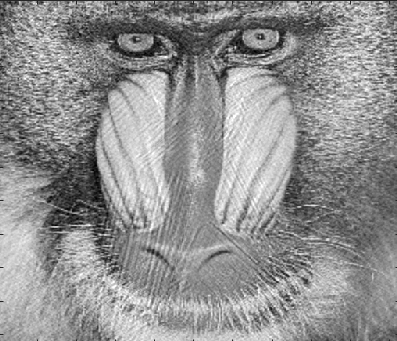}
    &
    \includegraphics[width=0.21\textwidth]{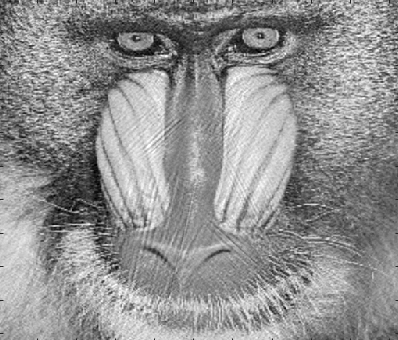}
    &
    \includegraphics[width=0.21\textwidth]{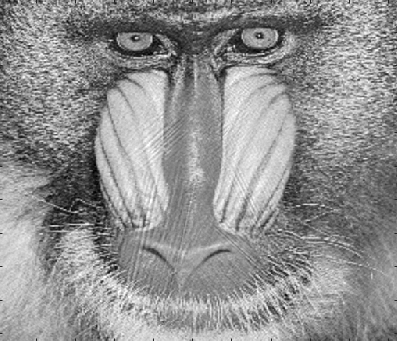}
    \\
    & 
    \rottext{phase}
    &
    \includegraphics[width=0.21\textwidth]{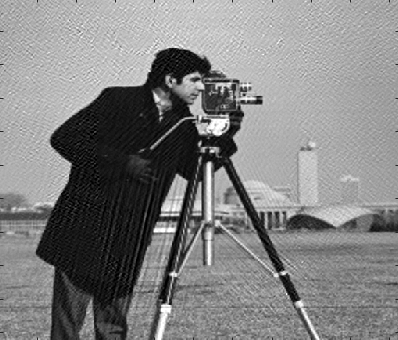}
    &
    \includegraphics[width=0.21\textwidth]{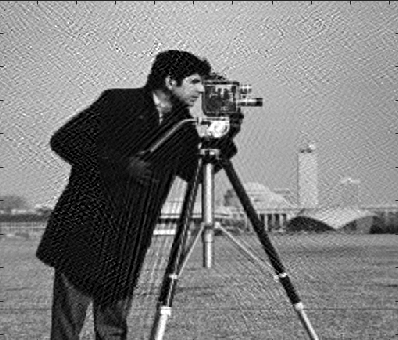}
    &
    \includegraphics[width=0.21\textwidth]{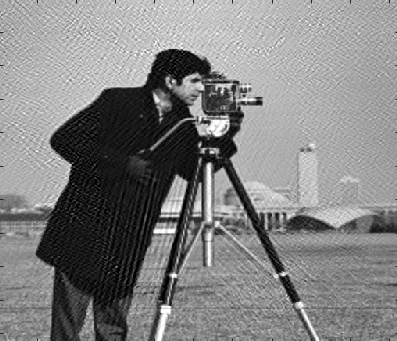}
    &
    \includegraphics[width=0.21\textwidth]{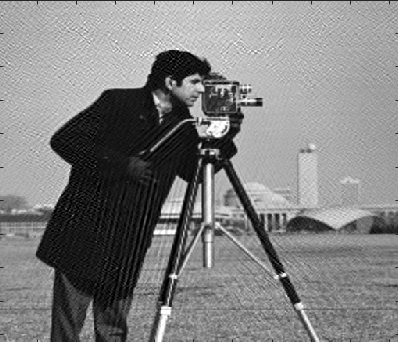}
    \\
    \hline
    \multirow{2}{*}{\rottextt{\textbf{$\mathbf{5\%}$ noise}}}
    &
    \rottext{magnitude}
    &
    \includegraphics[width=0.21\textwidth]{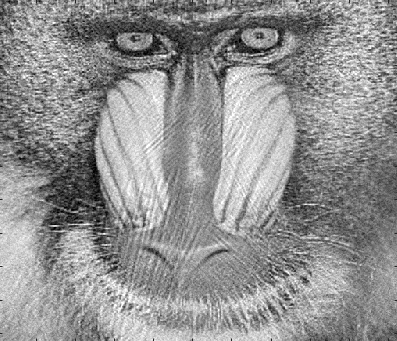}
    &
    \includegraphics[width=0.21\textwidth]{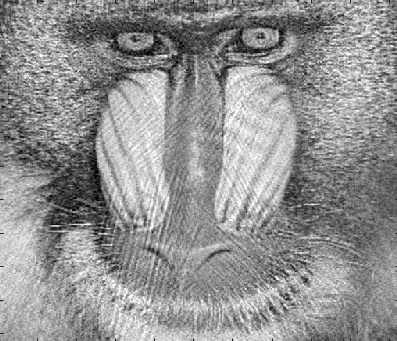}
    &
    \includegraphics[width=0.21\textwidth]{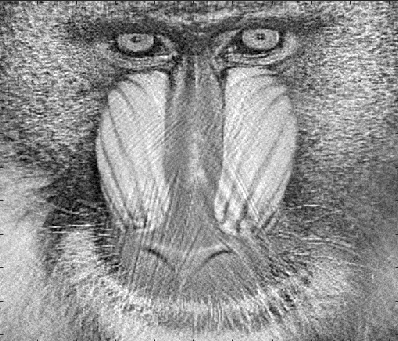}
    &
    \includegraphics[width=0.21\textwidth]{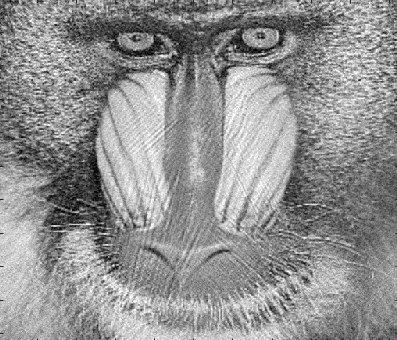}
    \\
    &
    \rottext{phase}
    &
    \includegraphics[width=0.21\textwidth]{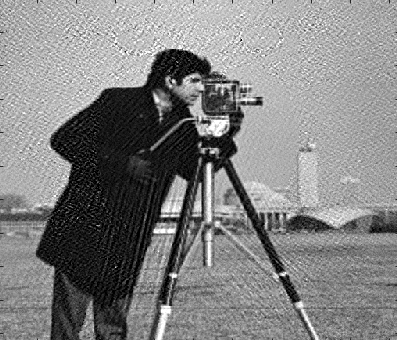}
    &
    \includegraphics[width=0.21\textwidth]{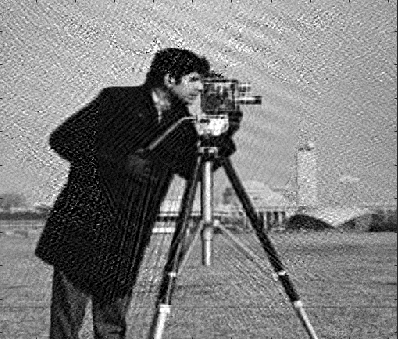}
    &
    \includegraphics[width=0.21\textwidth]{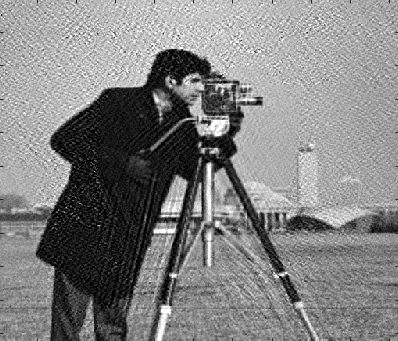}
    &
    \includegraphics[width=0.21\textwidth]{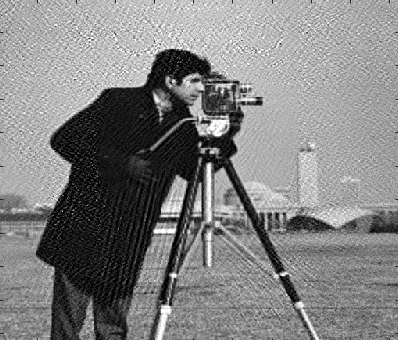}
    \\
    \hline
    \multirow{2}{*}{\rottextt{\textbf{$\mathbf{10\%}$ noise}}}
    &
    \rottext{magnitude}
    &
    \includegraphics[width=0.21\textwidth]{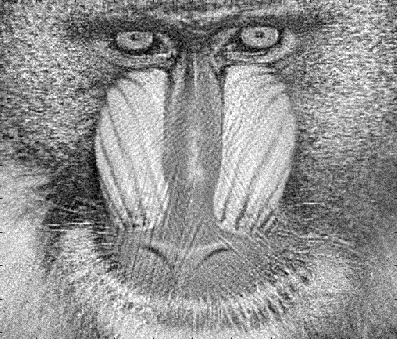}
    &
    \includegraphics[width=0.21\textwidth]{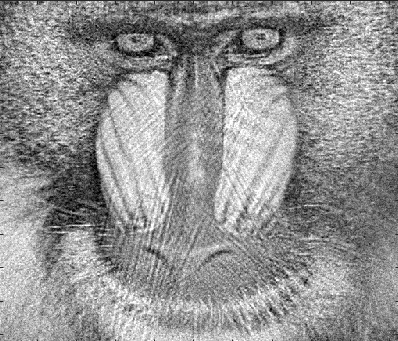}
    &
    \includegraphics[width=0.21\textwidth]{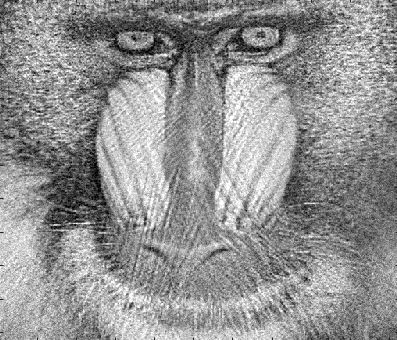}
    &
    \includegraphics[width=0.21\textwidth]{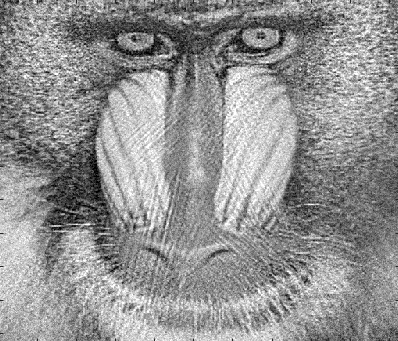}
    \\
    &
    \rottext{phase}
    &
    \includegraphics[width=0.21\textwidth]{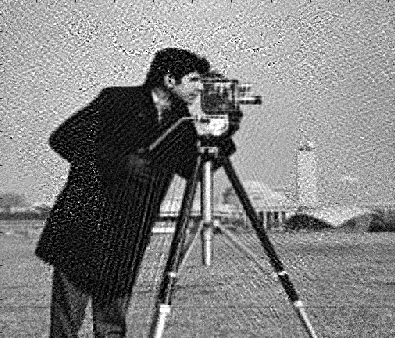}
    &
    \includegraphics[width=0.21\textwidth]{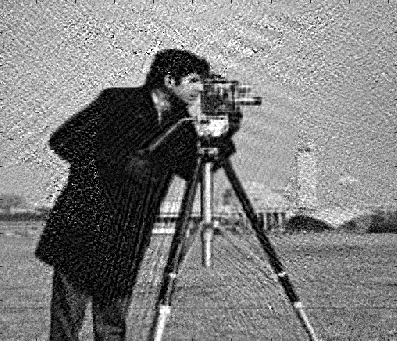}
    &
    \includegraphics[width=0.21\textwidth]{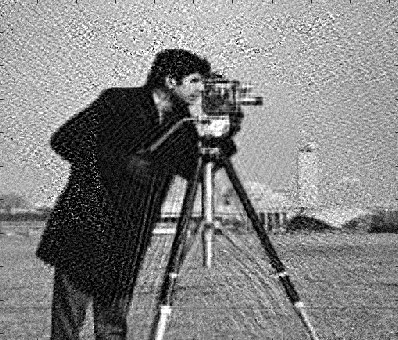}
    &
    \includegraphics[width=0.21\textwidth]{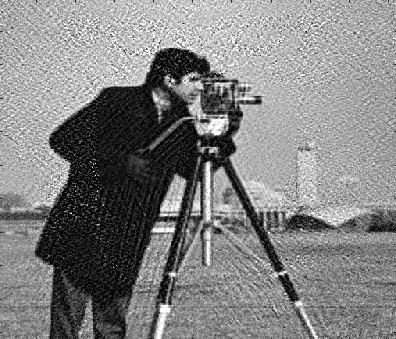}
    \end{tabular}
    \caption{\textit{Reconstructions of the magnitude (baboon) and phase (cameraman) for different noise levels.}}
    \label{fig:reconstructions}
    \vspace{-5mm}
\end{figure}
\section{Discussion}\label{sec:conclusion}
We present a multigrid approach for solving large-scale ptychographic phase retrieval problems. To this end, we consider the Multigrid Optimization (MG/OPT) scheme in order to improve the performance of gradient-based optimization algorithms for the ptychographic phase retrieval. The MG/OPT framework provides plenty of flexibility as it allows us to exploit the different hierarchical structures that ptychography exhibits.  Our approach is motivated by the Full Approximation Scheme, however, rather than attempting to solve a nonlinear system of equations, MG/OPT is a multigrid optimization framework where the coarse-grid subproblem is a first-order approximation to the fine grid problem; this guarantees a descent direction in the coarse-grid correction, and ultimately, convergence when a linesearch is performed (see Sec.~\ref{sec:MGOpt}) \cite{nash2010convergence}.

As expected, our numerical results show that MG/OPT improves the convergence of LBFGS for the ptychographic phase retrieval problem. Our numerical experiments also show that a 2-level and a 5-level MG/OPT with LBFGS as its underlying solver is competitive when compared to the Ptychographic Iterative Engine (PIE). In particular, the 5-level MG/OPT outperforms all other algorithms, which motivates more aggressive coarsening of the grids for large-scale problems. MG/OPT reduces the computational costs and accelerates the convergence of the ptychographic phase retrieval by shifting a substantial amount of work to the coarser grid. The acceleration of the convergence is particularly evident in the early iterations, when there is the most contribution coming from the coarser grids.
Further benefits are to be expected for larger problems, where deeper hierarchical structures can be exploited. We intend to extend our work for blind ptychography \cite{hesse2015proximal,maiden2009improved,nashed2014parallel,thibault2009probe}, where the probes must also be recovered, as well as for 3D ptycho-tomography \cite{gursoy2017direct,gilles20183d}.

\section*{Acknowledgments}
We thank Stefan Wild and Doga Gursoy for valuable discussions in the preparation of this paper. We also thank Meily Wu Fung for designing the schematic illustration in Fig.~\ref{fig:ptychoExperiment}. This material is supported by the U.S. Department of Energy, Office of Science, under contract DE-AC02-06CH11357.

\FloatBarrier
\bibliographystyle{abbrv}
\bibliography{references}



\end{document}